\documentclass[pdflatex,sn-mathphys-num]{sn-jnl}

\usepackage{graphicx}%
\usepackage{multirow}%
\usepackage{amsmath,amssymb,amsfonts}%
\usepackage{amsthm}%
\usepackage{mathrsfs}%
\usepackage[title]{appendix}%
\usepackage{xcolor}%
\usepackage{textcomp}%
\usepackage{manyfoot}%
\usepackage{booktabs}%
\usepackage{algorithm}%
\usepackage{algorithmicx}%
\usepackage{algpseudocode}%
\usepackage{listings}%
\usepackage{enumitem}
\usepackage{tcolorbox}
\usepackage{subcaption} 
\usepackage{multirow}
\usepackage[capitalise, nosort]{cleveref}
\crefname{equation}{}{}
\crefname{lemma}{Lemma}{Lemmas}
\crefname{assum}{Assumption}{Assumptions}
\crefname{corollary}{Corollary}{Corollaries}

\theoremstyle{thmstyleone}%
\newtheorem{theorem}{Theorem}[section]
\theoremstyle{thmstyletwo}%
\newtheorem{example}{Example}%
\newtheorem{remark}{Remark}%

\theoremstyle{thmstylethree}%
\newtheorem{definition}{Definition}[section]%

\newtheorem{assum}{Assumption}[section]
\newtheorem{lemma}{Lemma}[section]

\newcommand{\proj}[0]{\mathop{\bf Proj}}
\newcommand{\conv}[1]{{\bf co}\left\{ {#1} \right\}}
\newcommand{\R}{{\mathbb R}}
\newcommand{\argmin}[0]{ {\mathop{{\rm  argmin}}}}

\newcommand{\FL}[1]{{\mathcal F_{L}}\left( {#1} \right)}
\begin{document}

\title[Article Title]{Fast Convergence of Multiobjective Inertial Gradient Systems with Time Scaling}


\author[1]{\fnm{Yingdong}\sur{Yin}}\email{yydyyds@sina.com}



\affil*[1]{\orgdiv{National Center for Applied Mathematics in Chongqing}, 
\orgname{Chongqing Normal University}, 
\orgaddress{
\city{Chongqing}, 
\postcode{404100},
\country{China}}}




\abstract{In multiobjective optimization, inertial gradient systems accelerate convergence toward weakly Pareto optimal solutions. To achieve even faster convergence, we introduce a multiobjective inertial gradient system with time scaling (MITS), formulated as a second-order differential equation comprising an inertial term, asymptotically vanishing damping, and a time-scaled gradient term. We first establish the existence of solution trajectories for MITS. Through Lyapunov analysis, we show that with suitable parameters, the trajectory attains a convergence rate of \(O(1/t^{2}\beta(t))\) with respect to a merit function, where \(\beta(t)\) is a time-scaling function. Specifically, choosing \(\beta(t)=t^{p}\) for \(0\leq p<\alpha-3\) yields the rate \(O(1/t^{2+p})\), enabling arbitrarily fast sublinear convergence by tuning \(p\). We also prove that the trajectory converges to a weakly Pareto optimal solution. Furthermore, an implicit discretization of MITS leads to a multiobjective inertial proximal point method (MIPP), whose iterates share the \(O(1/k^{2}\beta_{k})\) rate and converge to a weakly Pareto optimum under appropriate conditions. Numerical experiments support the theoretical findings.}

\maketitle
\section{Introduction}
In this paper, $\mathbb{R}^n$ denotes the Euclidean space equipped with the inner product $\langle \cdot ,\cdot \rangle$ and the norm $\|\cdot \|$. Given a finite collection of convex functions $\{f_i:\mathbb{R}^n\to \mathbb{R}:i\in[m]:=\{1,2,\cdots,m\}\}$, we consider the following unconstrained multiobjective optimization problem:
\begin{equation} \label{eq:MOP}
\min_{x\in \mathbb{R}^n} F(x):=\left(f_1(x),f_2(x),\cdots,f_m(x)\right)^\top . \tag{MOP}
\end{equation}
Solving \cref{eq:MOP} in this paper refers to finding a weakly Pareto optimal solution, a point $x^*$ such that there exists no $y\in \mathbb{R}^n$ satisfying $f_i(y)<f_i(x^*)$ for all $i\in[m]$.

{\cref{eq:MOP} can be solved via multiobjective iterative methods. When the objective functions $\{f_i:i\in[m]\}$ are smooth, the core idea of these methods is to compute each new iterate using information from previous iterates together with a search direction-typically one that incorporates gradient information. {In this manner, they produce an infinite sequence whose cluster points correspond to weakly Pareto optimal solutions. Such methods originate from the multiobjective steepest descent method introduced by Fliege-Svaiter \cite{fliege2000steepest}. For a constant step-size $s$, the iteration is given by:}
\begin{equation}\label{eq:SD-scheme} 
x_{k+1}=x_k-s \proj_{C(x_k)}(0),
\end{equation}
where $C(x_k ) = \conv{\nabla f_i(x_k):i\in[m]}$ denotes the convex hull of the gradients of all objective functions, and $\proj_{C}(0)$ stands for the projection of the origin onto the closed convex set $C$. {For more multiobjective gradient-type iterative methods, refer to \cite{chen2023barzilai,chen2023conditional,assunccao2021conditional,fan2025faster,fliege2009newton,lucambio2018nonlinear,morovati2016barzilai,povalej2014quasi,qu2011quasi,tanabe2019proximal,tanabe2023accelerated,tanabe2022globally}.}

It is noteworthy that for smooth objective functions, we can adopt a dynamical systems perspective to study multiobjective gradient-type methods. In \cite{Attouch2014}, Attouch-Goudou pointed out that as the step-size $s$ tends to zero, \cref{eq:SD-scheme} can be transformed into a first-order gradient system as follows: 
\begin{equation}\label{eq:MOG} 
\dot x(t)+\proj_{C(x(t))}(0)= 0.\tag{MOG}
\end{equation}
The solution trajectory of \cref{eq:MOG} and the iterative sequence generated by \cref{eq:SD-scheme} share similar convergence properties. Under the characterization by the merit function proposed in \cite{tanabe2024new}, they achieve convergence rates of $O(1/t)$ and $O(1/k)$, respectively, as shown in \cite{sonntagphdthesis,YingdongYin2026BalancedGradientFlow}. To achieve a faster convergence rate, based on the inertial gradient system proposed by Attouch-Garrigos \cite{attouch2015multiibjective}, Sonntag-Peitz \cite{sonntag2024fastSIAMonOptimization} introduced the following \textit{Multiobjective Inertial Gradient System with Asymptotic Vanishing Damping}:
\begin{equation} \label{eq:MAVD} 
\frac{\alpha }{t}\dot x(t)+\proj_{C(x(t))+\ddot x(t)}(0) = 0.\tag{MAVD}
\end{equation}
For $\alpha \ge 3$, under the same merit function characterization, the solution trajectory of \cref{eq:MAVD} achieves a convergence rate of $O(1/t^2)$ and converges to a weakly Pareto optimal solution \cite{yin2025pointConvergenceAnalysis}. Further studies on multiobjective gradient systems can be found in \cite{boct2024inertial,luo2025accelerated,YingdongYin2026BalancedGradientFlow}. For $m=1$, \cref{eq:MAVD} reduces to the following inertial gradient system proposed by Su-Boyd-Cand\`es \cite{su2016differential} for single-objective optimization:
\begin{equation}\label{eq:single-optimization-gradient-system}
\ddot x(t)+\frac{\alpha }{t}\dot x(t)+\nabla f_1(x(t)) = 0.
\end{equation}
As in the multiobjective setting, for $\alpha \ge 3$, the solution of \cref{eq:single-optimization-gradient-system} possesses a convergence rate of $O(1/t^2)$ and converges to the optimal solution \cite{Jang2025Point,boct2025iteratesnesterovsacceleratedalgorithm,attouch2018fast}.

I n the single-objective optimization setting, time scaling can lead to faster convergence of the gradient system's trajectory. Attouch-Chbani-Riahi \cite{attouch2019fast,attouch2019fast0} proposed the following \textit{Inertial Gradient System with Time Scaling}:
\begin{equation}\label{eq:ITS}
\ddot x(t)+\frac{\alpha }{t}\dot x(t)+\beta(t)\nabla f_1(x(t))=0,\tag{ITS}
\end{equation} 
where $\alpha\ge 3$ and $\beta :\mathbb{R}\to \mathbb{R}$ is an increasing smooth function with $\lim _{t\to \infty }\beta(t) = +\infty $. When $\dot \beta(t)\le (\alpha -3)\beta(t)/t$, the convergence rate of the \cref{eq:ITS} trajectory in terms of function values is $O(1/t^2\beta(t))$. In particular, a convergence rate of $O(1/t^{2+p})$ can be achieved when $\beta(t)=t^p$ with $0\le p<\alpha -3$. 
More on dynamical systems with time scaling of the gradient can be found in \cite{adly2024accelerated,attouch2025fastclosedloop,attouch2019fast0,he2022secondorderprimalfirstorderdual,attouchandbot2024fast,csetnek2024second}. Based on an implicit discretization of \cref{eq:ITS}, for the case where $f_1$ is lower semicontinuous, \cite{attouch2019fast} presents an inertial (accelerated) proximal point method:
\begin{equation}\label{eq:iter-1}
\left\{\begin{aligned}
	y_k &=x_k+\frac{k-1}{k+\alpha -1}(x_k-x_{k-1}),\\
	x_{k+1}&=\textbf{prox}_{\lambda_k f_1}(y_k):= \argmin_{z\in\mathbb{R}^n}\left\{f_1(z)+\frac{1}{2\lambda_{k}}\|z-y_k\|^2\right\},
\end{aligned}\right.
\end{equation}
where $\alpha \ge 1$, and
\begin{equation} 
\lambda _k = \frac{k\beta_k }{k+\alpha -1} ,\qquad \beta_{k+1}\le \frac{k(k+\alpha -1)}{(k+1)^2}\beta_k.
\end{equation}
This algorithm achieves a function value convergence rate of $O(1/k^2\beta_k)$. Specifically, for $\alpha > 3$ and $\beta_k = k^p$ with $0 \le p < \alpha - 3$, the algorithm achieves an $O(1/k^{2+p})$ convergence rate, matching the trajectory convergence rate in \cref{eq:ITS}.

The above observations motivate us to investigate analogous properties for gradient systems and proximal point methods in the case $m>1$, i.e., in multiobjective optimization. Specifically, we study the following two questions:
\begin{itemize}
	\item \textbf{Question-1}: Is it possible to apply time scaling to the gradient for \cref{eq:MIPP} accelerate the convergence rate of the merit function along the solution trajectory?
	\item \textbf{Question-2}: Is it possible to construct a multiobjective proximal point method similar to the iteration scheme \cref{eq:iter-1}, such that the convergence rate of the merit function along the generated sequence is faster than $O(1/t^2)$?
\end{itemize}
\subsection{Main contribution for Question-1}
By $\proj_{C+x}(0)=\proj_C(-x)+x$, we consider an equivalent form of \cref{eq:MAVD}:
\begin{equation}\label{eq:MAVD-other} 
\ddot x(t)+\frac{\alpha }{t}\dot x(t)+\proj_{C(x(t))}(-\ddot x(t)) = 0.
\end{equation} 
Suppose that the twice differentiable function $X(t)$ satisfies \cref{eq:MAVD-other}. If we apply the time scaling $t = s^p$ and define $Y(s) := X(s^p)$, we obtain:
\begin{equation}\label{eq:time-rescale}
\dot Y(s)=ps^{p-1}\dot X(s^p),\qquad \ddot Y(s)=p(p-1)s^{p-1}\dot X(s^p)+p^2 s^{2(p-1)}\ddot X(s^p).
\end{equation}
Based on \cref{eq:MAVD-other,eq:time-rescale}, we deduce that $Y(s)$ is a solution to the following equation:
\begin{equation}\label{eq:system-example}
\ddot y(s)+\frac{\alpha_p}{s}\dot y(s)+{\proj_{p^2s^{2(p-1)}C(y(s))}\left(-\ddot y(s)+\frac{p-1}{s}\dot y(s)\right)}=0,
\end{equation}
where $\alpha_p :=1+(\alpha - 1 )p$. Since $\|\dot X(t)\| = O(1/t)$, according to the first equality in \cref{eq:time-rescale}, we obtain $\|\dot Y(s)\| = O(1/s)$. This observation inspires us to construct a simpler system. Firstly, consider the following function:
\begin{equation*} 
\varphi(s): = \proj_{p^2s^{2(p-1)}C(Y(s))}\left(-\ddot Y(s)\right)-\proj_{p^2s^{2(p-1)}C(Y(s))}\left(-\ddot Y(s)+\frac{p-1}{s}\dot Y(s)\right).
\end{equation*}
Moreover, by the non-expansiveness of the projection mapping, we obtain
\begin{equation*} 
\|\varphi(s)\|\le \frac{p-1}{s}\|\dot Y(s)\|=O(1/s^2). 
\end{equation*} 
Therefore, based on \cref{eq:time-rescale} and the specific form of $\varphi(s)$, and by once again using the identity $\proj_{C}(-x) + x = \proj_{C + x}(0)$,  $Y(s)$ satisfies the following equation:
\begin{equation}\label{eq:MITS-yesterday}
\frac{\alpha _p}{s}\dot y(s)+\proj_{p^2 s^{2(p-1)}C(y(s))+\ddot y(s)}(0) = \varphi(s)=o(1).
\end{equation}
Based on \cref{eq:MITS-yesterday} and neglecting the perturbations on the right-hand side, we propose the following {\it Multiobjective Inertial Gradient System with Time Scaling}:
\begin{equation}\label{eq:MITS}
\frac{\alpha }{t}\dot x(s)+\proj_{\beta(t)C(x(t))+\ddot x(t)}(0) = 0.\tag{MITS}
\end{equation}
For \cref{eq:MITS}, we obtain the following main results:
\begin{itemize}
    \item \textbf{Existence of solution trajectories}. We prove the existence of solution trajectories for the \cref{eq:MITS} within a general framework. Specifically, there exists a smooth function $x(t)$ defined on $[t_0, +\infty)$ that satisfies the given initial conditions, and for almost all $t \in [t_0, +\infty)$, $x(t)$ satisfies \cref{eq:MITS}. 
    \item \textbf{Arbitrarily fast convergence rate}. Under certain assumptions on $\beta(t)$ and utilizing the merit function, we derive the convergence rate of the function values along the solution trajectory is $O(1/t^2\beta(t))$. In particular, $\beta(t) = t^p$ can satisfy the assumptions, allowing the convergence rate to vary with $p$ and achieve arbitrarily fast rates within the range of sublinear convergence.
    \item \textbf{Convergence of solution trajectories}. We prove that the solution trajectory $x(t)$ converges to a weakly Pareto optimal solution $x^*$ of \cref{eq:MOP}. To the best of our knowledge, this is the latest result for inertial gradient systems with time scaling, including the single-objective optimization case.
    \item \textbf{Faster convergence rate with a new metric}. Using $x^*$, we show that $\mathcal W(t):=\min_{i\in[m]}(f_i(x(t))-f_i(x^*))+\frac{1}{2 t^p}\|\dot x(t)\|^2$ achieves a convergence rate of $o(1/t^{2+p})$. This measure characterizes, to some extent, the convergence rate of the objective function values along the trajectory toward the optimal value.
\end{itemize}
\subsection{Main contribution for Question-2}
For any selected sequence $\{t_k\}$ with $t_k \to \infty$ as $k\to \infty$, \cref{eq:MITS} implies the following:  
\begin{equation}\label{eq:main-question-2-1}
\ddot x(t_k)+\frac{\alpha  - 1}{t_k}\dot x(t_k)+\frac{1}{t_k}\dot x(t_k)+\beta(t_k)\cdot \conv{\nabla f_i(x(t_k)):i\in [m]}\ni0.
\end{equation} 
Setting $t_k = k$, $\beta_k :=\beta(t_k)$, $x_{k+1}:=x(t_k)$, and  
\begin{equation} \label{eq:main-question-2-2}
\ddot x(t_k )\approx x_{k+1}-2x_k+x_{k-1}, \qquad \dot x(t_k )\approx x_{k+1}-x_k \approx x_{k}-x_{k-1}.
\end{equation} 
Under the assumption that $\{f_i:i\in[m]\}$ are lower semi-continuous but not smooth, \cref{eq:main-question-2-2} motivates us to construct the iterative sequence $\{x_k\}$ of a proximal point method via the following relation:  
\begin{equation*} 
\frac{1}{\beta_{k}}\bigg[x_{k+1}-2x_k+x_{k-1}+\frac{\alpha -1}{k}(x_{k+1}-x_k)+\frac{1}{k}(x_k-x_{k-1})\bigg] + \conv{\bigcup_{i\in[m]}\partial f_i(x_{k+1})}\ni 0,
\end{equation*} 
where $\partial f_i(x)$ denotes the subdifferential set of $f_i$ at $x$ for $i \in [m]$. After simplification, the above can be equivalently written as  
\begin{equation}\label{eq:iter-yesterday} 
\frac{1}{\lambda _k }\left(x_{k+1}-y_k \right)+\conv{\bigcup_{i\in[m]}\partial f_i(x_{k+1})}\ni0,  
\end{equation} 
where $y_k = x_{k}+\frac{k-1}{k+\alpha -1}(x_{k}-x_{k-1})$ and $\lambda_k  = \frac{k\beta_k}{k+\alpha -1}$. When $m = 1$, \cref{eq:iter-yesterday} corresponds to the first-order optimality condition for the subproblem in \cref{eq:iter-1}. Moreover, for any $\xi = (\xi_1, \xi_2, \cdots, \xi_m)^\top \in \mathbb{R}^m$, the inclusion  
$\partial \max_{i \in [m]} (f_i(z) - \xi_i) \subseteq \conv{\bigcup_{i \in [m]} \partial f_i(z)}$
holds for all $z \in \mathbb{R}^n$ (\cite[Theorem 3.50]{beck2017first}), which motivates us to consider setting up a minimax subproblem to construct the desired proximal method. In fact, Tanabe-Fukuda-Yamashita constructed a class of multiobjective proximal methods in \cite{tanabe2023accelerated}:  
\begin{equation}\label{eq:Tanabe-Fukuda-Yamashita-method} 
\left\{\begin{aligned}
y_k &=x_k+\gamma_k (x_k-x_{k-1}),\\
x_{k+1}& =\underset{z\in \mathbb{R}^n}{\argmin}\left\{\max_{i\in[m]}(f_i(z)-f_i(x_k))+\frac1{2s}\|z-y_k\|^2\right\},
\end{aligned}\right.
\end{equation} 
where $\gamma_k$ is a relevant parameter. It is easy to verify that when $\gamma_k =\alpha_k$ and $\lambda_k=s$, the iterates generated by the above satisfy \cref{eq:iter-yesterday}. Based on \cref{eq:Tanabe-Fukuda-Yamashita-method}, we define $\Phi_{x_k}(z):=\max_{i\in[m]}(f_i(z)-f_i(x_k))$ and construct the following \textit{Multiobjective Inertial Proximal Point Method}:  
\begin{equation}\label{eq:MIPP} 
\left\{\begin{aligned}
y_k &=x_k+\frac{k-1}{k+\alpha -1}(x_k-x_{k-1}),~~\lambda _k =\frac{k\beta_k}{k+\alpha-1}, \\
x_{k+1}&=\textbf{prox}_{\lambda_k \Phi_{x_k}}(y_k)=\underset{z\in \mathbb{R}^n}{\argmin}\left\{\Phi_{x_k}(z)+\frac1{2\lambda_k }\|z-y_k\|^2\right\},
\end{aligned}\right.\tag{MIPP}
\end{equation}
where $\{\beta_k\}$ is a given increasing positive sequence. Clearly, when $m=1$, \cref{eq:MIPP} is reduced to the iterative scheme \cref{eq:iter-1}. For \cref{eq:MIPP}, our main contributions are as follows: 
\begin{itemize}
    \item \textbf{Arbitrarily fast convergence rate of the iterations}. For the generated sequence $\{x_k\}$ of \cref{eq:MIPP}, we obtain a convergence rate of $O(1/k^2\beta_k)$ using the merit function. {Numerical experiments reflect this result.}
    \item \textbf{Convergence of iterations}. {Under certain conditions, we prove that $\{x_k\}$ converges to a weakly Pareto optimal solution of \cref{eq:MOP}.} 
\end{itemize}
\subsection{Organization}
The organization of this paper is as follows: In \cref{sec:pre}, we provide the necessary preliminaries; \cref{sec:exist} contains the main results on the existence of \cref{eq:MITS} solution trajectories and their basic properties; in \cref{sec:asymptotic}, we analyze the convergence rate and convergence of the solution trajectory; \cref{sec:algo} presents the convergence rate and convergence analysis for \cref{eq:MIPP}; \cref{sec:num} carries out numerical experiments to validate the results related to \cref{eq:MIPP}; \cref{sec:con} offers conclusions and future prospects; the Appendix provides a more detailed discussion of the existence of solution trajectories and the corresponding proofs.
\section{Preliminaries}\label{sec:pre}
\subsection{Notation}
In this paper, $\mathbb{R}^d$ denotes a $d$-dimensional Euclidean space with the scalar product $\langle \cdot, \cdot \rangle$ and the induced norm $\|\cdot\|$. 
{For any positive integer $m$, let $[m]:=\{1,2,\cdots,m\}$.} 
For any vectors $a, b \in \mathbb{R}^d$, we say $a \leq b$ if $a_i \leq b_i$ holds for all $i\in[d]$; the relations $<$, $\geq$, and $>$ are defined analogously. The nonnegative  and positive orthants are denoted respectively by
$\mathbb{R}_{+}^d:=\{x\in \mathbb R^d : x\geq 0\}, ~~\mathbb{R}^d_{++}:=\{x\in \mathbb{R}^d:x>0\}.$
The open and closed balls of radius $\delta>0$ centered at a point $x\in \mathbb{R}^d$ are defined, respectively, as:
$\mathcal  B_{\delta}(x) := \left\{ y \in \mathbb{R}^d : \|y - x\| < \delta \right\}, ~~\overline{\mathcal  B}_{\delta}(x) := \left\{ y \in \mathbb{R}^d : \|y - x\| \le \delta \right\}. $
The set $\Delta^m := \{ \theta \in \mathbb{R}^m : \theta \geq 0 \text{ and } \sum_{i=1}^m \theta_i = 1 \}$ is the positive unit simplex. Given a set of vectors $\{\eta_1, \ldots, \eta_m\} \subseteq \mathbb{R}^d$, their convex hull is defined as $\conv {\eta_1, \ldots, \eta_m} := \{ \sum_{i=1}^m \theta_i \eta_i : \theta \in \Delta^m \}$. For a closed convex set $C \subseteq \mathbb{R}^d$, the projection of a vector $x$ onto $C$ is  $\proj_{ C}(x) := \argmin_{y \in C} \|y - x\|^2$. 
Given a vector-valued function $F: \mathbb{R}^n \to \mathbb{R}^m$, the lower level set of $F$ for a given level $a \in \mathbb{R}^m$ is defined as   
$\mathcal L(F,a):=\left\{x\in \mathbb{R}^n: F(x)\leq a\right\}.$ For a constant $ L > 0 $, the set $\FL{\mathbb{R}^d}$ is referred to as the class of $L$-smooth functions on $\mathbb{R}^d$, i.e. for any function $ f \in \FL{\mathbb{R}^d} $ defined on $\mathbb{R}^d$, we have  
\begin{equation*} 
\left\|\nabla f(x) - \nabla f(y)\right\| \le L \|x - y\|, \quad \forall x, y \in \mathbb{R}^d.
\end{equation*}
\subsection{Pareto optimality}
The concepts and  optimal condition for \cref{eq:MOP} are detailed as follows:
\begin{definition}\label{def:defofpareto}
Consider the multiobjective optimization problem \cref{eq:MOP}.
\begin{enumerate}[label=(\roman*)]        
\item A point $ x^* \in \mathbb{R}^n $ is called a weakly Pareto optimal or weakly Pareto optimal solution if there has no $y\in \mathbb R^n$ that $F(y)<F(x^*)$. The set of all weakly Pareto optimal solutions is called the weak Pareto set and is denoted by $ \mathcal{P}_w $. The image $F(\mathcal P_w)$ of the Pareto set $\mathcal P_w$ is the weak Pareto front. A point $x^*\in \mathbb{R}^n$ is a locally Pareto optimal solution if there exists $\delta >0$ such that $x^*$ is Pareto optimal in $\mathcal B_{\delta}(x^*)$. 
\item A point $ x^* \in \mathbb{R}^n $ is called Pareto optimal or a Pareto optimal solution  if there has no $y\in \mathbb R^n$ that $F(y)\leq F(x^*)$ and $F(y)\neq F(x^*)$. The set of all Pareto optimal solutions is called the Pareto set and is denoted by $ \mathcal{P} $.  The image $F(\mathcal P)$ of the Pareto set $\mathcal P$ is the Pareto front. A point $x^*\in \mathbb{R}^n$ is a locally weakly Pareto optimal solution if there exists $\delta >0$ such that $x^*$ is weakly Pareto optimal in $\mathcal B_{\delta}(x^*)$. 
\end{enumerate}  
\end{definition}
\begin{definition}\label{def:KKTcondition}
    Assume that $f_i$ is continuously differentiable for $i\in[m]$.  A point $ x^* \in \mathbb{R}^n $ is said to satisfy the Karush-Kuhn-Tucker (KKT) conditions if
    \begin{equation}\label{eq:KKTpoint}
    \proj_{\conv{\nabla f_i(x(t)):i\in[m]}}(0) = 0 .
    \end{equation}
    If $ x^* $ satisfies the KKT conditions, it is called Pareto critical or a Pareto critical point. The set of all Pareto critical points is called the Pareto critical set and is denoted by $ \mathcal{P}_c $.
\end{definition}	
\begin{lemma}
    The following statements hold:
    \begin{enumerate}[label=(\roman*)]
        \item If $x^*\in \mathbb{R}^n$ is locally weakly Pareto optimal for \cref{eq:MOP}, then $x^*\in\mathcal P_c$;
        \item If each component function $f_i$ ($i\in[m]$) is convex and $x^*\in \mathcal P_c$, then $x^*\in \mathcal P_w$ ;
        \item If $f$ is strictly convex and $x^*$ is Pareto critical for \cref{eq:MOP}, then $x^*\in \mathcal P$. 
    \end{enumerate}
\end{lemma}
\subsection{Merit function}
A merit function is defined as a nonnegative function in \cref{eq:MOP} that vanishes only at weakly Pareto optimal solutions. In this work, we also examine the merit function introduced by Tanabe-Fukuda-Yamashita. \cite{tanabe2024new}:
\begin{equation}\label{eq:meritfunction} 
	u_0(x) := \sup_{z \in \mathbb{R}^n} \min_{i\in [m]} \bigl( f_i(x) - f_i(z) \bigr). 
\end{equation}
This function plays a role analogous to $ f(x) - f(x^*) $ in single-objective optimization-it is nonnegative and serves as an indicator of weak optimality, as rigorously stated in the following theorem.
\begin{theorem}[{\cite[Theorem 3.1, 3.2]{tanabe2024new}}]\label{thm:weakpareto} 
    Let $ u_0:\mathbb R^n \to \mathbb{R}\cup{\{\pm \infty\} } $ be defined as in \cref{eq:meritfunction}. Then,
    \begin{enumerate}[label=(\roman*)] 
        \item $ u_0(x) \geq 0 $ for all $ x \in \mathbb{R}^n $;
        \item  $ x \in \mathbb{R}^n $ is a weakly Pareto optimal solution of \cref{eq:MOP} if and only if $ u_0(x) = 0 $;
        \item $ u_0(x) $ is lower semi-continuous.
\end{enumerate}
\end{theorem}
For any $a\in \R^n$, we denote ${\mathcal L}{\mathcal P}_w(F,a):=\mathcal P_w\cap {\mathcal L}(F,a)$.  
\begin{theorem}[{\cite[Lemma 2.7]{sonntag2024fastSIAMonOptimization}}]\label{thm:sup-inf}
Let $x_0\in \mathbb{R}^n$ and $x\in\mathcal L(F,F(x_0))$, then 
\begin{equation*}
u_0(x) =	\sup_{F^*\in F( \mathcal L\mathcal P_w(F,F(x_0)))} \inf_{z\in F^{-1}(F^*)} \min_{i\in[m]} (f_i(x)-f_i(z)).
\end{equation*}
\end{theorem}
\begin{assum}\label{assume:boundedsup}
    For any $ x_0 \in \mathbb{R}^n $ and all $ x \in \mathcal{L}(F, F(x_0)) $, there exists $ x^* \in \mathcal{L}\mathcal{P}_w (F, F(x_0)) $ such that  
    \begin{equation}
    \mathbf{R} := \sup_{F^* \in F\left(\mathcal{L}\mathcal{P}_w(F, F(x_0))\right)} \inf_{z \in F^{-1}(F^*)} \frac12\|z - x_0\|^2 < +\infty.
    \end{equation}
\end{assum} 
\section{Existence of solutions and auxiliary lemmas}\label{sec:exist}
In this section, we  define the solution trajectory of \cref{eq:MITS} and present auxiliary lemmas useful for asymptotic analysis.
\subsection{Existence of solution trajectory}
 Similar to the setup in \cite{boct2024inertial}, in this paper, smooth functions satisfying \cref{eq:MITS} that possess a series of properties are referred to as solution trajectory, with the specific definition given as follows:

\begin{definition} \label{def:soultion}
    A function $x:[t_0,+\infty )\to \mathbb{R}^n$, $t\mapsto x(t) $, is called a solution or a solution trajectory to \cref{eq:MITS} if it satisfies the following conditions:  
    \begin{enumerate}[label=(\roman*)] 
        \item $x(\cdot)\in C^1([t_0,+\infty ))$.  
        \item $\dot x(t)$ is absolutely continuous on $[t_0,T]$ for any $T\ge t_0$.  
        \item There exists a measurable function $\ddot x(t):[t_0,+\infty )\to \mathbb{R}^n$ such that $\dot x(t)=\dot x(t_0)+\int_{t_0}^t\ddot x(s)ds$ for $t\ge t_0$.  
        \item $\dot x(\cdot) $ is almost everywhere differentiable and $\frac{d}{dt}\dot x(t)=\ddot x(t)$ holds for almost all $t\in [t_0,+\infty )$.  
        \item $\ddot x(t)+\frac{\alpha }{t}\dot x(t)+{\proj}_{\beta(t) C(x(t))}(-\ddot x(t))=0$ for almost all $t\in[t_0,+\infty )$.
        \item $x(t_0)=x_0$ and $\dot x(t_0)=v_0$.  
    \end{enumerate}
\end{definition}
We provide a detailed discussion of solution existence in the \cref{appendix:existence}. Specifically, for the \cref{eq:MITS}, we have the following result:
\begin{theorem} 
     Suppose that $f_i\in \FL{\R^n}$ for all $i\in[m]$, then there exists a function $x:[t_0,+\infty )\to \R^n$ which is a solution of \cref{eq:MITS} in the sense of \cref{def:soultion}
\end{theorem}   
\subsection{Auxiliary functions and lemma}  
The following auxiliary lemmas can be regarded as intrinsic properties of the solution trajectory $ x(t) $ itself, and they are useful in asymptotic analysis.
\begin{itemize}
    \item  For each $i\in [m]$, define the energy function:  
    \begin{equation}\label{eq:aux-fun-wi} 
    \mathcal W_i:[1,+\infty )\to \mathbb{R},\quad t\mapsto f_i(x(t))+\frac{1}{2\beta(t) }\|\dot x(t)\|^2.
    \end{equation}  
    \item For $z\in \mathbb{R}^n$, define the auxiliary function
    \begin{equation}\label{eq:aux-fun-Theta}
    \Theta_z(t):[1,+\infty )\to \mathbb{R},~~~t\mapsto \min_{i\in[m]} \left(f_i(x(t))-f_i(z)\right)
    \end{equation}
\end{itemize}
The following are the fundamental properties of the solution trajectory:
\begin{lemma}\label{lem:energy}
Suppose that $f_i\in \FL{\R^n}$ for all $i\in[m]$. Let $x:[1,+\infty )\to\R^n$ be a solution to \cref{eq:MITS} with $(x(1),\dot x(1)) = (x_0,0)$. Denote $\Omega_{\alpha}(t):=\frac{1}{\beta(t)}\left(\frac{\alpha}{t}+\frac{\dot \beta(t)}{2}\right)$. Then, we have
    \begin{enumerate}[label=(\roman*)]
        \item $\mathcal W_i(t)$ is non-increasing and 
        \begin{equation}\label{eq:W-diff-0}\frac{d}{dt}\mathcal W_i(t)+\Omega_{\alpha}(t)\|\dot x(t)\|^2 \le 0, \end{equation} for almost all $t\in[1,+\infty )$ and $i\in[m]$. 
        \item The following inequality holds for any $t''\ge t' \ge 1$:
        \begin{equation}
        -\Theta_z(t') \le -\Theta_z(t'')-\frac{1}{2\beta(t'')}\|\dot x(t'')\|^2 +\frac{1}{2\beta(t')}\|\dot x(t')\|^2-\int_{t'}^{t''} \Omega_{\alpha}(\tau )\|\dot x(\tau)\|^2d\tau.
        \end{equation}
        \item $x(t)\in \mathcal L(F,F(x_0))$ for all $t\in [1,+\infty )$. 
        \item Assuming $f_i$ is bounded below, then $\lim_{t\to \infty }\mathcal W_i(t)=\mathcal W^*\in \mathbb{R}$ exists and  
        \begin{equation*} 
        \int_{1}^{+\infty} \Omega_{\alpha}(t)\|\dot x(t)\|^2dt<+\infty.
        \end{equation*}
    \end{enumerate} 
\end{lemma}  
\begin{proof} 
(i) 
Since $\mathcal{W}_i(t)$ is almost everywhere differentiable on $[1, +\infty)$, for almost all $t \in [1, +\infty)$, we have  
\begin{equation}\label{eq:W-diff}
\begin{aligned}
\frac{d}{dt}\mathcal W_i(t)&=\left<\nabla f_i(x(t)),\dot x(t)\right>+\frac{1}{\beta(t)}\left\langle\ddot x(t),\dot x(t)\right\rangle-\frac{\dot \beta(t)}{2\beta(t)}\|\dot x(t)\|^2\\
&=\frac{1}{\beta(t)}\left[\left\langle \beta(t)\nabla f_i(x(t))+\ddot x(t),\dot x(t)\right\rangle \right] -\frac{\dot \beta(t)}{2\beta(t)}\|\dot x(t)\|^2.
\end{aligned}
\end{equation} 
By definition of \cref{eq:MITS} and the projection theorem, we have  
\begin{equation}\label{eq:W-projection} 
\left\langle\beta(t)\nabla f_i(x(t))+\ddot x(t),\dot x(t) \right\rangle\le -\frac{\alpha }{t}\left\|\dot x(t)\right\|^2.
\end{equation}  
Combining \cref{eq:W-diff} and \cref{eq:W-projection}, we obtain the result. 

(ii) According to (i), for any $1 \le t_1 \le t_2$,  
\begin{equation*} 
\begin{aligned}
\mathcal{W}_i(t_2) - \mathcal{W}_i(t_1) &= \int_{t_1}^{t_2} \frac{d}{d\tau} \mathcal{W}_i(\tau) \, d\tau \le -{\int_{t_1}^{t_2} \Omega(\tau)\|\dot x(\tau)\|^2 \, d\tau}=:{\mathbf{I}(t_1,t_2)},
\end{aligned}
\end{equation*}
holds for every $i \in [m]$.  Moreover, from the definition of $\mathcal{W}_i$, taking any $z \in \mathbb{R}^n$, we obtain by direct computation  
\begin{equation*} 
-\bigl(f_i(x(t_1)) - f_i(z)\bigr) \le -\bigl(f_i(x(t_2)) - f_i(z)\bigr) - \frac{1}{2\beta(t_2)}\|\dot x(t_2)\|^2 + \frac{1}{2\beta(t_1)}\|\dot x(t_1)\|^2 - \mathbf{I}(t_1,t_2),
\end{equation*}
which holds for every $i \in [m]$. Therefore, the conclusion follows.

(iii) Based on (i) and the fact that $\dot x(1) = 0$, the conclusion follows.

(iv) 
Since $\mathcal{W}_i(t)$ is bounded from below and from (i), it follows that $\lim_{t \to \infty} \mathcal{W}_i(t) = \inf_{t \geq 1} \mathcal{W}_i(t) := \mathcal{W}_i^*$ exists. Based on this result and combining with \cref{eq:W-diff-0}, we have  
\begin{equation*} 
\int_{1}^{+\infty} \frac{1}{\beta(t)}\left(\frac\alpha t+\frac{\dot \beta(t)}{2}\right) dt \leq \mathcal{W}_i(1) - \mathcal{W}_i^* < +\infty.
\end{equation*} 
Therefore, the conclusion holds.
\end{proof}
\section{Asymptotic and convergence analysis}\label{sec:asymptotic}
For a continuously differentiable function $\beta:[1,+\infty)\to \mathbb{R}$, define the following functions  
\begin{equation}\label{eq:Omega-Gamma}
\Omega_{\alpha}(t) := \frac{1}{\beta(t)} \left( \frac{\alpha}{t} + \frac{\dot{\beta}(t)}{2} \right), \quad \Gamma_{\alpha}(t) := (\alpha - 3) t \beta(t) - t^2 \dot{\beta}(t).
\end{equation} 
Consider the following assumptions:
\begin{assum}\label{assum:beta-t}
{Assume that the smooth function $\beta(t)$ is positive and non-decreasing on $[1,+\infty)$}  
\begin{enumerate}[label=(\roman*)]
\item { There exists $t_1 \ge 1$ such that $\Gamma_{\alpha}(t) \ge 0$ for all $t \ge t_1$;}  
\item {There exist $t_2 \ge t_1$ and $\eta > 0$ such that the following inequality holds for all $t \ge t_2$: } 
\begin{equation*} 
 \frac{\Gamma_{\alpha}(t)}{2\beta(t)} - \Omega_{\alpha}(t) \int_{t_1}^{t} \Gamma_{\alpha}(s) \, ds \le -\eta \cdot \frac{\Gamma_{\alpha}(t)}{2\beta(t)}.
\end{equation*}   
\end{enumerate}
\end{assum}
\begin{example}\label{example:beta}
Let the function $\beta(t) = \beta_0 t^p (\ln t)^{q(\alpha -3-p)}$, where $\beta_0 > 0$, $0 \le p \le \alpha - 3$, $q \ge 0$; and $r_{\alpha ,p,q}(t):= \alpha - p - 3 - \frac{q(\alpha -3-p)}{\ln t}$. We have
\begin{equation}\label{eq:exmaple-Gamma_alpha}
\Gamma_{\alpha}(t) = \beta_0 t^{p+1} (\ln t)^{q(\alpha-3-p)} \, r_{\alpha,p,q}(t), \quad 
\frac{\Gamma_{\alpha}(t)}{2\beta(t)} = \frac{t}{2} \, r_{\alpha ,p,q}(t),
\end{equation}
\begin{equation}\label{eq:example-Omega_alpha}
\Omega_{\alpha}(t) = \frac{\alpha}{\beta_0 t^{p+1} (\ln t)^{q(\alpha -3 - p)}} + \frac{p}{2t} + \frac{q(\alpha -3-p)}{2t\ln t}. 
\end{equation}
Hence there exists $t_1 \ge 1$ such that for any $t \ge t_1$, $\Gamma_{\alpha}(t) \ge 0$, i.e., $\beta(t)$ satisfies \cref{assum:beta-t}(i). It is easy to compute that there exists $\bar t \ge t_1$ such that for any $t \in [\bar t, +\infty)$,  
\begin{equation*}
	\int_{t_1}^t \Gamma_{\alpha}(s) \, ds \ge \frac23 \, \frac{\beta_0}{p+2} \, t^{p+2} (\ln t)^{q(\alpha - 3 - p)} \, r_{\alpha ,p,q}(t).
\end{equation*}
Combining with \cref{eq:example-Omega_alpha}, we obtain  
\begin{equation}\label{eq:example-r}
\begin{aligned}
	\Omega_{\alpha}(t) \int_{t_1}^{t} \Gamma_{\alpha}(s) ds 
	&\ge \frac{t}{p+2} \left[ \frac23 \alpha + \frac{\beta_0 p}{3} \, t^p (\ln t)^{q(\alpha - 3 - p)} \right] r_{\alpha,p,q}(t).
\end{aligned}
\end{equation}
When $p = 0$, since $\alpha \ge 3$, \cref{eq:example-r} implies that
\begin{equation}\label{eq:example-Gamma-beta}
\Omega_{\alpha}(t) \int_{t_1}^t \Gamma_{\alpha}(s) \, ds \ge 2 \cdot \frac{\Gamma_{\alpha}(t)}{2\beta(t)},
\end{equation}
for all $t \in [\bar t, +\infty)$. When $p > 0$, there exists $t_2 \ge \bar t$ such that \cref{eq:example-Gamma-beta} holds for all $t\ge t_2$. Therefore, $\beta(t)$ satisfies \cref{assum:beta-t}(ii).
\end{example}
\subsection{Fast convergence}\label{subsec:fastconvergence}
Constructing auxiliary energy functions is a general method for deriving the convergence rates of solution trajectories. For single-objective gradient systems, such functions are nonnegative and referred to as Lyapunov functions. However, in the context of multiobjective optimization, this property does not hold \cite{luo2025accelerated,sonntag2024fastSIAMonOptimization,boct2025fast}. Therefore, we refer to the constructed function as a Lyapunov-like function, as follows:
\begin{equation}\label{eq:aux-fun-lya-like}
\mathcal{E}_z(t) = \underbrace{t^{2}\beta(t)\Theta_z(t)}_{\textbf{Potential energy }\mathcal E_{z}^{\rm pot}(t)} + \underbrace{\frac{1}{2}\|(\alpha -1 )(x(t) - z) + t\dot{x}(t)\|^2}_{\textbf{Mixed energy }\mathcal E_{z}^{\rm mix}(t)}.
\end{equation}
The following lemma shows that the derivative of $\mathcal E_z(t)$ is dominated by $-\Theta_z(t)$.
\begin{lemma}\label{lem:non-increase-of-Lyapunov}
Suppose that $f_i\in \FL{\R^n}$ for all $i\in[m]$. Let $x:[t_0,+\infty )\to \mathbb{R}^n$ be a solution to \cref{eq:MITS} with $(x(t),\dot x(t))=(x_0,0)$. Let $\mathcal{E}_z(t)$ be defined as in \cref{eq:aux-fun-lya-like}, then 
\begin{equation*}
\frac{d}{dt}\mathcal E_z(t)\le-\Gamma_{\alpha }(t)\Theta_z(t),
\end{equation*}
for almost all $t\in[1,+\infty )$, where $\Gamma_{\alpha}(t)$ is defined in \cref{eq:Omega-Gamma}. 
\end{lemma}
\begin{proof}
According to the definition of \cref{eq:MITS}, there exists $(\theta_1(t),\cdots,\theta_m(t))\in \Delta^m$ such that  
\begin{equation}
\frac{\alpha }{t}\dot x(t)+\ddot x(t)=-\proj_{\beta(t) C(x(t))}(-\ddot x(t))=-\beta(t)\sum_{i=1}^{m}\theta_i(t)\nabla f_i(x(t)).
\end{equation}
According to the definition of \cref{eq:MITS} again and by direct computation, we have  
\begin{equation}\label{eq:aux-fun-lya-like-diff-1}
\begin{aligned}
\frac{d}{dt}\mathcal E_{z}^{\rm mix}(t)&  = t^2\left<\dot x(t),\ddot x(t)+\frac{\alpha}{t}\dot x(t)\right>+(\alpha -1)t\left<x(t)-z,-\beta(t)\sum_{i=1}^{m}\theta_i(t)\nabla f_i(x(t))\right>\\
&\le t^2\left<\dot x(t),\ddot x(t)+\frac{\alpha}{t}\dot x(t)\right>-(\alpha -1)t\beta(t)\Theta_z(t), 
\end{aligned}
\end{equation} 
{where the inequality holds for the convexity of $f_i$.} Since $\Theta_z(t)$ is differentiable for almost all $t\in[1,+\infty )$, and according to \cite[Lemma 4.12]{sonntag2024fastSIAMonOptimization}, there exists $i_t\in[m]$ such that  
\begin{equation}\label{eq:aux-fun-lya-like-diff-2}
\frac{d}{dt}\mathcal E_z^{\rm pot}(t)= (2t\beta(t)+t^2\dot \beta(t))\Theta_z(t)+t^{2}\Big<\dot x(t),\beta(t) \nabla f_{i_t}(x(t))\Big>.
\end{equation}
Combining \cref{eq:aux-fun-lya-like-diff-1} and \cref{eq:aux-fun-lya-like-diff-2}, we obtain for almost all $t\in[1,+\infty )$,  
\begin{equation*} 
\frac{d}{dt}\mathcal E_z(t)+t((\alpha -3)\beta(t)-t\dot \beta(t))\Theta_z(t)\le t^2\Big<\dot x(t),\frac{\alpha }{t}\dot x(t)+\ddot x(t)+\beta(t) \nabla f_{i_t}(x(t))\Big>\le 0,
\end{equation*}  
where the last inequality follows from the projection theorem.
\end{proof}
\begin{lemma}\label{lem:Lya-inequal-1}
 Suppose that $f_i\in \FL{\R^n}$ for all $i\in[m]$, and \cref{assum:beta-t} holds true. Let $x:[t_0,+\infty )\to \mathbb{R}^n$ be a solution to \cref{eq:MITS} with $(x(t),\dot x(t))=(x_0,0)$. Let $\mathcal{E}_z(t)$ be defined as in \cref{eq:aux-fun-lya-like}, then
\begin{equation*}
\mathcal E_z(t) +\Theta_z(t) \int_{t_1}^t\Gamma_{\alpha}(s)ds+\int_{t_1}^t\frac{\Gamma_{\alpha}(s)}{2\beta(t)}\|\dot x(s)\|^2ds+\eta\int_{t_1}^t \|\dot x(\tau)\|^2d\tau \le \mathcal E_z(1)+\mathbf{I}_3+\mathcal J_z,
\end{equation*}
holds for all $t\in[t_2,+\infty )$, where $\mathbf{I_3} := \int_{t_1}^{t_2}\left|\frac{\Gamma(\tau)}{2\beta(\tau)}-\Omega(\tau)\int_{t_1}^{\tau}\Gamma(s)ds\right| \|\dot x(\tau)\|^2 d\tau$ and  $\mathcal J_z:= \int_{1}^{t_1}|\Gamma_{\alpha}(s)|\cdot |\Theta_z(s)|ds $. 
\end{lemma}
\begin{proof}
To simplify notation, we denote $\Gamma$ and $\Omega$ instead of $\Gamma_{\alpha}$ and $\Omega_{\alpha}$. According to \cref{lem:non-increase-of-Lyapunov}, we have  
\begin{equation}
\label{eq:lya-1}\frac{d}{ds}\mathcal E_z(s)\le -\Gamma(s)\Theta_z(s),
\end{equation}
for all $s\in [1,+\infty )$. {According to \cref{lem:energy} and $\Gamma(s)\ge 0$ for all $s\geq t_1$, we have},
\begin{equation*}
\begin{aligned}
-\Gamma(s)\Theta_z(s) \le -\Gamma(s)\Theta_z(t)-\frac{\Gamma(s)}{2\beta(t)}\|\dot x(t)\|^2+\frac{\Gamma(s)}{2\beta(s)}\|\dot x(s)\|^2-\Gamma(s)\int_{s}^{t}\Omega (\tau)\|\dot x(\tau)\|^2d\tau, 
\end{aligned}
\end{equation*}
holds for any $t\ge s\ge t_1$. Moreover, integrating over $[1,t]$ for \cref{eq:lya-1},  
\begin{equation}\label{eq:lya-2}
\begin{aligned}
	-\int_{1}^{t}\Gamma(s)\Theta_z(s)ds
	& \le -\int_{1}^{t_1}\Gamma(s)\Theta_z(s)ds-\Theta_z(t)\int_{t_1}^t\Gamma(s)ds-\frac{\int_{t_1}^t\Gamma(s)ds}{2\beta(t)}\|\dot x(t)\|^2\\
	&\qquad + \underbrace{\int_{t_1}^t \frac{\Gamma(s)}{2\beta(s)}\|\dot x(s)\|^2ds}_{\mathbf{I}_1(t)} -\underbrace{\int_{t_1}^{t}\Gamma(s)\int_{s}^t\Omega(\tau)\|\dot x(\tau)\|^2d\tau ds}_{\mathbf{I}_2(t)}. 
\end{aligned}
\end{equation}  
By Fubini's theorem, interchanging the order of integration in $\mathbf{I}_2(t)$ yields  
\begin{equation*} 
\begin{aligned}
\mathbf{I}_2(t)&=\int_{t_1}^t\Omega(\tau)\|\dot x(\tau)\|^2\int_{t_1}^\tau\Gamma_{\alpha}(s)\ dsd\tau,
\end{aligned}
\end{equation*}  
and hence  
\begin{equation}\label{eq:lya-3}
\begin{aligned}
	\mathbf{I}_1(t)  -\mathbf{I}_2(t) &= \left(\int_{t_1}^{t_2}+\int_{t_2}^t\right)\left[\frac{\Gamma(\tau)}{2\beta(\tau)}-\Omega(\tau)\int_{t_1}^{\tau}\Gamma(s)ds\right]
	\|\dot x(\tau)\|^2d\tau\\
	& \le \int_{t_1}^{t_2}\left|\frac{\Gamma(\tau)}{2\beta(\tau)}-\Omega(\tau)\int_{t_1}^{\tau}\Gamma(s)ds\right|\cdot \|\dot x(\tau)\|^2 d\tau - \eta\int_{t_2}^{t}\frac{\Gamma(\tau)}{2\beta(\tau)}\|\dot x(\tau)\|^2d\tau.
\end{aligned}
\end{equation}  
Therefore, integrating over $[1,t]$ for \cref{eq:lya-1} and combining with \cref{eq:lya-2,eq:lya-3}, we obtain  
\begin{equation*} 
\mathcal E_z(t) -\mathcal E_z(1)+\Theta_z(t)\int_{t_1}^{t} \Gamma(s)ds+\eta \int_{t_1}^{t} \frac{\Gamma(\tau)}{2\beta(\tau)}\|\dot x(\tau)\|^2d\tau
\le \mathbf{I}_3+\mathcal{J}_z.
\end{equation*}
We complete the proof. 
\end{proof}
Based on the lemmas above, we present the convergence rate.
\begin{theorem}\label{thm:fastcpnvergence}
    Suppose that $f_i\in \FL{\R^n}$ for all $i\in[m]$, and \cref{assume:boundedsup,assum:beta-t} hold true. Let $\alpha \geq 3$, and let $x: [1, +\infty) \to \mathbb{R}^n$ be a bounded solution to \cref{eq:MITS} with $(x(1),\dot x(1))=(x_0,0)$. Then,
\end{theorem}
\begin{enumerate}[label=(\roman*)]
    \item  $u_0(x(t))=O(1/t^{2}\beta(t))$.
    \item  $\|\dot x(t)\|=O(1/t)$.
    \item  $\frac{\Gamma_{\alpha}(t)}{2\beta(t)} \|\dot x(t)\|^2\in L([1,+\infty ),\mathbb{R})$.
\end{enumerate}
\begin{proof}
According to  \cref{lem:Lya-inequal-1},  we obtain  
\begin{equation}\label{eq:Theta_z-1}
\begin{aligned}
&~~~~\mathcal E_{z}(t)+
\Theta_z(t)\underbrace{ \left(t^2\beta(t)+\int_{t_1}^t\Gamma_{\alpha}(s)ds\right)}_{\mathbf{I}_4(t)}+\underbrace{\eta\int_{t_1}^t\Gamma_{\alpha}(\tau)\|\dot x(\tau)\|^2d\tau}_{\mathbf{I}_5(t)}\\&\le \mathcal E_z(1)+\mathbf{I}_3+\mathcal{J}_z\le u_0(x_0)+\frac{(\alpha -1)^2}2\|x_0-z\|^2 +\mathbf{I}_3+\mathcal{J}_z,
\end{aligned}
\end{equation} 
for all $z\in\mathbb{R}^n$. Furthermore, based on the convexity of \(f_i\), we obtain  
\begin{equation}\label{eq:J_z}
\begin{aligned}
	\mathcal J_z
	&\leq \int_{1}^{t_1}|\Gamma_{\alpha}(s)|\cdot \min_{i\in[m]}\|\nabla f_i(x(s))\|\|x(s)-z\|ds\\
	&\leq \int_{1}^{t_1}|\Gamma_{\alpha}(s)|\cdot \min_{i\in[m]}\|\nabla f_i(x(s))\|\Big(\|x(s)-x_0\|\Big)ds\\
	&\qquad+\frac{1}{2}\left(\int_{1}^{t_1}|\Gamma_{\alpha}(s)|\cdot \min_{i\in[m]}\|\nabla f_i(x(s))\|ds\right)^2 +\frac12\|x_0-z\|^2\\
	&\leq \underbrace{\mathbf{M}\cdot \left(\sup_{t\ge 1}\|x(t)\|+\|x_0\|\right)+\frac12{\mathbf{M}^2}}_{\widetilde {\mathbf{M}}}+\frac12\|x_0-z\|^2,
\end{aligned}
\end{equation}
where $\mathbf{M}= \int_{1}^{t_1} |\Gamma_{\alpha}(s)|\cdot \min_{i\in[m]}\|\nabla f_i(x(s))\| ds$. Thus, applying \cref{thm:sup-inf}, we obtain  
\begin{equation}\label{eq:Theta_z-2}
\begin{aligned} 
t^2\beta(t) u_0(x(t))+\mathbf I_{5}(t)&\le\mathbf{I}_4(t) \sup_{F^*\in \mathcal L\mathcal P_w(F,F(x_0))}\inf_{z\in F^{-1}(F^*)}\Theta_z(t)+\mathbf{I}_5(t)\\
&\leq u_0(x_0) +[(\alpha -1)^2+1]\mathbf{R}+\widetilde{\mathbf{M}}+\mathbf{I}_3,
\end{aligned}
\end{equation}
and hence $u_0(x(t))=O(1/t^2\beta(t))$. 
Moreover, let $z= x(t)$ for \cref{eq:Theta_z-1,eq:J_z}, we have
\begin{equation}\label{eq:Theta_z-3}
{\frac{1}{2}t^2\|\dot x(t)\|^2\le u_0(x_0) +[(\alpha -1)^2+1] \left( \sup_{t \ge 1}\|x(t)\|^2 + \|x_0\|^2 \right)+\widetilde {\mathbf{M}} +\mathbf{I}_3=: \widetilde {\mathbf R},} 
\end{equation}
i.e. $\|\dot x(t)\|=O(1/t)$. 
Since $u_0(x(t))\ge 0$, according to \cref{eq:Theta_z-2},  we have
\begin{equation*}
\eta \int_{\bar t}^{+\infty}\frac{\Gamma_{\alpha}(t)}{2\beta(t)}dt=\lim_{t\to +\infty}\mathbf{I}_5(t) <+\infty. 
\end{equation*}
This completes the proof.
\end{proof}
\subsection{Trajectory convergence analysis}\label{subsec:convergenceofsoultion}
The convergence of solution trajectories for \cref{eq:ITS,eq:MITS} under specific parameter selections remained an open problem for a long time and was recently resolved \cite{Jang2025Point,yin2025pointConvergenceAnalysis}. Building on recent studies, we provide a proof of convergence for the solution trajectories of \cref{eq:MITS}.
\begin{lemma}\label{lem:continuous-Lyapunov-limit-exist}
 Suppose that $f_i\in \FL{\R^n}$ for all $i\in[m]$, \cref{assume:boundedsup,assum:beta-t} hold true. Let  $x: [1, +\infty) \to \mathbb{R}^n$ be a bounded solution to \cref{eq:MITS} with $(x(1),\dot x(1))=(x_0,0)$. Assuming $f_i$ is bounded below for $i\in[m]$, then 
\begin{enumerate}[label=(\roman*)]
    \item The limit $\lim_{t \to \infty} f_i(x(t)) = f_i^\infty$ exists;
    \item For any cluster point $z^*$ of $x(t)$, we have
    \begin{equation*} 
   \lim_{t \to \infty} \mathcal E_{z^*}(t) = \mathcal E^\infty_{z^*},
    \end{equation*}
    exists.
\end{enumerate}
\end{lemma}
\begin{proof}
 (i) Based on the non-increasing property of $\mathcal W_i(t)$, it is deduced that for any $t \ge s \ge 1$,   
 \begin{equation*}  
 \mathcal W(t) := u_0(x(t)) + \frac{1}{2\beta(t)}\|\dot x(t)\|^2 \le u_0(x(s)) + \frac{1}{2\beta(s)}\|\dot x(t)\|^2.
 \end{equation*}    
 Thus, $\mathcal W(t)$ is non-increasing. Since $u_0(x(t)) \ge 0$, the limit $\lim_{t \to \infty} \mathcal W(t) = \inf_{t\ge 1} \mathcal W(t) := \mathcal W^*$ exists. Furthermore, as $u_0(x(t)) \to 0$ as $t \to \infty$ by \cref{thm:fastcpnvergence}, we have   
 \begin{equation*}  
 \lim_{t \to \infty} \frac{1}{2\beta(t)}\|\dot x(t)\|^2 = \mathcal W^* \in \mathbb{R}. 
 \end{equation*}     
 By \cref{lem:energy} (ii), it follows that $\lim_{t \to \infty} \frac{1}{2\beta(t)}\|\dot x(t)\|^2 = 0$. Consequently,   
 \begin{equation*}
\lim_{t \to \infty} f_i(x(t)) = \lim_{t \to \infty} \left( \mathcal W_i(t) - \frac{1}{2\beta(t)}\|\dot x(t)\|^2 \right) = \mathcal W_i^\infty =: f_i^\infty \in \mathbb{R},
 \end{equation*}   
 exists.   
 
(ii) From (i), for  any cluster point $z^*$ of $x(t)$, we obtain   
\begin{equation} \label{eq:bouded-Theta_z}
\begin{aligned} \Theta_{z^*}(t) &= \min_{i\in[m]} \left( \mathcal W_i(t) - \frac{1}{2\beta(t)}\|\dot x(t)\|^2 - f_i(z^*) \right) \ge -\frac{1}{2\beta(t)}\|\dot x(t)\|^2. \end{aligned} 
\end{equation}   
Based on \cref{thm:fastcpnvergence} and \cref{eq:Theta_z-3}, it follows that 
\begin{equation*}
{t^{2}\beta(t)\Theta_{z^*}(t) \ge t^{2} \cdot \left( -\frac{1}{2}\|\dot x(t)\|^2 \right) \ge -\widetilde {\mathbf R}.}
\end{equation*}   
Thus, $\mathcal E_{z^*}(t) \ge -\widetilde {\mathbf R}$, i.e., it is bounded below. Moreover, according to \cref{lem:non-increase-of-Lyapunov}, we have
\begin{equation*} 
{\frac{d}{dt}\mathcal E_{z^*}(t)\le \Gamma_{\alpha }(t)\cdot \frac{1}{2\beta(t)}\|\dot x(t)\|^2,}
\end{equation*}
where the right-hand side of the inequality is integrable on $[1,+\infty)$ by \cref{thm:fastcpnvergence}.
Consequently, $\lim_{t \to \infty} \mathcal E_{z^*}(t) = \mathcal E_{z^*}^\infty$ exists.  
\end{proof}
{Some words in here.}
\begin{lemma}[{\cite[Lemma A.4]{boct2025fast}}]\label{lem:continuous-limit} Let $a>0$ and $q:[t_0,+\infty )\to \R^n$ be a continuously differentiable function such that 
    \begin{equation}
    \lim_{t \to \infty}\left(q(t)+\frac{t}{a}\dot q(t)\right) = \ell\in \R.
    \end{equation}
    Then, it holds $\lim_{t\to \infty}q(t) = \ell$. 
\end{lemma}
\begin{theorem}\label{thm:point-convergence-continuous}
 Suppose that $f_i\in \FL{\R^n}$ for all $i\in[m]$, \cref{assume:boundedsup} hold true. Let  $x: [1, +\infty) \to \mathbb{R}^n$ be a bounded solution to \cref{eq:MITS} with $(x(1),\dot x(1))=(x_0,0)$. Assuming $f_i$ is bounded from below for $i\in[m]$, then $\lim_{t\to \infty} x(t) = x^*\in \mathbb{R}^n$ exists and $x^*$ is a weakly Pareto optimal solution of \cref{eq:MOP}.  
\end{theorem}
\begin{proof}
For any $ z \in \mathbb{R}^n $, define
\begin{equation*} 
h_z(t) = \frac12\|x(t)-z\|^2.
\end{equation*}
Let $ z_1 $ and $ z_2 $ be two cluster points of $ x(t) $. According to the definition of $ \mathcal E_z(t) $, we have
\begin{equation*} 
\begin{aligned}
\mathcal{E}_{z_j}(t)
&=t^2\beta(t)\Theta_{z_j}(t)+\frac{t^2}2\|\dot x(t)\|^2+\frac{(\alpha -1)^2}2\|x(t)-z_j\|^2+(\alpha-1)t\left\langle \dot x(t),x(t)-z_j \right\rangle \\
&=t^2\beta(t)\Theta_{z_j}(t)+\frac{t^2}2\|\dot x(t)\|^2+(\alpha -1)^2h_{z_j}(t) +(\alpha -1) t \dot h_{z_j}(t),\quad j = 1, 2.
\end{aligned}
\end{equation*} 
Thus,
\begin{equation*} 
\begin{aligned}
\frac{1}{(\alpha -1)^2} \left( \mathcal E_{z_1}(t) - \mathcal E_{z_2}(t) \right) &= \left[ h_{z_1}(t) - h_{z_2}(t) \right] + \frac{t}{\alpha -1} \left[ \dot h_{z_1}(t) - \dot h_{z_2}(t) \right].
\end{aligned}
\end{equation*}  
By \cref{lem:continuous-Lyapunov-limit-exist}, since $ \mathcal E_{z_j}(t) $ has a limit for $j = 1, 2$, the left-hand side of the above equation has a limit, denoted as $ \ell $. According to \cref{lem:continuous-limit}, we have
\begin{equation*} 
\frac12 \lim_{t\to \infty } \left[ \|x(t)-z_1\|^2 - \|x(t)-z_2\|^2 \right] = \lim_{t\to \infty } \left[ h_{z_1}(t) - h_{z_2}(t) \right] = \ell.
\end{equation*}
Take sequences $ \{t_k\} $ and $ \{s_k\} $ such that $ x(t_k) \to z_1 $ and $ x(s_k) \to z_2 $ as $ k \to \infty $. Then
\begin{equation*} 
\begin{aligned}
\lim_{k\to \infty } \left[ \|x(t_k)-z_1\|^2 - \|x(t_k)-z_2\|^2 \right] &= -\|z_1 - z_2\|^2 = 2\ell, \\
\lim_{k\to \infty } \left[ \|x(s_k)-z_1\|^2 - \|x(s_k)-z_2\|^2 \right] &= \|z_1 - z_2\|^2 = 2\ell,
\end{aligned}
\end{equation*}
so $ z_1 = z_2 := x^* $, i.e., $ \lim_{t\to \infty } x(t) = x^* $. By \cref{thm:weakpareto} (iii), we obtain
\begin{equation*} 
0 \le u_0(x^*) \le \liminf_{t\to \infty } u_0(x(t)) = 0,
\end{equation*}
which implies $ u_0(x^*) = 0 $. According to \cref{thm:weakpareto} (ii), $ x^* $ is a weakly Pareto optimal solution of \cref{eq:MOP}.
\end{proof}
\subsection{Faster convergence with a gap function}
In this subsection, we choose $\beta(t) = t^p$ with $p < \alpha - 3$ to present an indicator that decreases faster than $O(1/t^{2+p})$. Under this setting, according to \cref{thm:fastcpnvergence,example:beta}, the following holds:
\begin{equation*} 
\int_{1}^{+\infty} t\|\dot x(t)\|^2 dt < +\infty.
\end{equation*} 
By \cref{thm:fastcpnvergence}, there exists a weakly Pareto optimal solution $x^*$ of \cref{eq:MOP} such that $\lim_{t\to \infty} x(t) = x^*$. Define
\begin{equation}\label{eq:energy-gap} 
{\mathcal{W}_{x^*}(t) := \min_{i \in [m]} \left( \mathcal{W}_i(t) - \lim_{t \to \infty} \mathcal{W}_i(t) \right) = \min_{i \in [m]} (f_i(x(t)) - f_i(x^*)) + \frac{1}{2t^p} \|\dot x(t)\|^2.}
\end{equation} 
By \cref{lem:energy}, it can be seen that $\mathcal{W}_{x^*}(t) \ge 0$ and is non-increasing. According to \cref{thm:fastcpnvergence}, we obtain
\begin{equation*}
\mathcal{W}_{x^*}(t) \le u_0(x(t)) + \frac{1}{2t^{p}} \|\dot x(t)\|^2 = O(1/t^{2+p}).
\end{equation*} 
The following theorem shows that the above convergence rate is not optimal.

\begin{theorem} \label{thm:small-o}
 Suppose that $f_i\in \FL{\R^n}$ and bounded from below for all $i\in[m]$, \cref{assume:boundedsup} hold true. Let  $x: [1, +\infty) \to \mathbb{R}^n$ be a bounded solution to \cref{eq:MITS} with $(x(1),\dot x(1))=(x_0,0)$, $\beta(t) = t^p$, $p\le \alpha-3$. Then
\begin{equation*} 
\mathcal{W}_{x^*}(t) = o\left( \frac{1}{t^{2+p}} \right).
\end{equation*} 
\end{theorem}
\begin{proof} Note that in this case, $\Gamma_{\alpha}(t) = (\alpha - 3 - p)t^{p+1}$. Furthermore, by \cref{lem:non-increase-of-Lyapunov}, we have
\begin{equation*} 
\int_{1}^{+\infty} \Gamma_{\alpha}(t) \Theta_{x^*}(t) dt \le -\int_{1}^{+\infty} \frac{d}{dt} \mathcal{E}_{x^*}(t) dt = \mathcal{E}_{x^*}(1) - \mathcal{E}^\infty < +\infty.
\end{equation*} 
Thus,
\begin{equation*} 
\int_1^{+\infty} \Gamma_{\alpha}(t) \mathcal{W}_{x^*}(t) dt = \int_{1}^{+\infty} \Gamma_{\alpha}(t) \Theta_{x^*}(t) dt + \frac{\xi}{2} \int_{1}^{+\infty} t \|\dot x(t)\|^2 dt < +\infty.
\end{equation*}
In addition, note that
\begin{equation*} 
\begin{aligned}
\mathcal{E}_{x^*}(t) &= t^{2+p} \Theta_{x^*}(t) + \frac{t^2}{2} \|\dot x(t)\|^2 + \frac{(\alpha - 1)^2}{2} \|x(t) - x^*\|^2 + (\alpha - 1) t \langle \dot x(t), x(t) - x^* \rangle \\
&= t^{2+p} \mathcal{W}_{x^*}(t) + \frac{(\alpha - 1)^2}{2} \|x(t) - x^*\|^2 + (\alpha - 1) t \langle \dot x(t), x(t) - x^* \rangle .
\end{aligned}
\end{equation*} 
Since $\|\dot x(t)\| = O(1/t)$ and $\lim_{t \to \infty} \|x(t) - x^*\| = 0$, it follows that $t \langle \dot x(t), x(t) - x^* \rangle \to 0$ as $t \to \infty$. Therefore,
\begin{equation*} 
\lim_{t \to \infty} t^{2+p} \mathcal{W}_{x^*}(t) = \lim_{t \to \infty} \mathcal{E}_{x^*}(t) = \mathcal{E}^\infty \in \mathbb{R}.
\end{equation*} 
Hence, $\lim_{t \to \infty} t^{2+p} \mathcal{W}_{x^*}(t) = 0$.
\end{proof}
\begin{remark}
For the case $ m = 1, p = 0 $, \cref{thm:small-o} yields the following result:  
\begin{equation}\label{eq:for-scaler-o}
\mathcal{W}_{x^*}(t) = f_1(x) - \inf_{x \in \mathbb{R}^n} f_1(x) + \frac{1}{2} \|\dot{x}(t)\|^2 = o\left(\frac{1}{t^2}\right).
\end{equation} 
Based on this conclusion, we obtain the classical results $ f_1(x(t)) - f_1(x^*) = o(1/t^2) $ and $ \|\dot{x}(t)\| = o(1/t) $ \cite{may2017asymptoticfasterrate,attouch2016rate}. However, when $ m > 1 $, we cannot derive similar conclusions from \cref{eq:for-scaler-o}, as the difficulty lies in the fact that $ \min_{i \in [m]} (f_i(x) - f_i(x^*)) $ is not necessarily non-negative. Indeed, if we already have the result $ \|\dot{x}(t)\| = o(1/t) $, then using \cref{thm:fastcpnvergence} and \cref{eq:for-scaler-o}, we can deduce $ \left| \min_{i \in [m]} (f_i(x(t)) - f_i(x^*)) \right| = o(1/t^2) $. Therefore, under the current assumptions, proving $ \|\dot{x}(t)\| = o(1/t) $ becomes critical; alternatively, constructing relevant counterexamples to negate this conclusion is an open problem we leave for future work.
\end{remark}
 
\section{Inertial proximal point method}\label{sec:algo}
{ In this section, we investigate the generated sequences of \cref{eq:MIPP} and present the corresponding convergence properties and rates under certain assumption.} In particular, we give an assumption for $\{\beta_k\}$: 
\begin{equation}\label{eq:beta_k}
\beta_{k+1}  = \frac{k(k+\alpha -1)}{(k+1)^2}\beta_k. \tag{$H_{\beta}$}
\end{equation}
Based on the above, we propose the following algorithm to generate the iterative sequence of \cref{eq:MIPP}:
\begin{algorithm}[H]
    \caption{Multiobjective Inertial Proximal Point Method}
    \label{algo:MIPP}
    \begin{algorithmic}[1] 
        \Require Initial values: $\alpha\geq 3$, $k_{\max}>0$,  $x_0=x_1\in \mathbb{R}^n$,  and $\beta_1>0$.
        \State Set $k=1$.
        \While{$k<k_{\max}$}
        \State Compute $\lambda_k$, $y_k$ and $x_{k+1}$ by \cref{eq:MIPP}. 
        \If{ $\|\frac{1}{\lambda_k}(x_{k+1}-y_k)\|<\varepsilon$}
        \State \Return $x_{k+1}$
        \Else
        \State Compute $\beta_{k+1}$ by \cref{eq:beta_k}. 
        \State Set $k\to k+1$
        \EndIf
        \EndWhile		
    \end{algorithmic}
\end{algorithm}
The sequence generated by \cref{algo:MIPP} or \cref{eq:MIPP}, as discussed in the subsequent text, refers to an infinite sequence without considering termination conditions. {The following lemma can be regarded as a fundamental property of the iterative sequence generated by \cref{algo:MIPP}.}
\begin{lemma}\label{lem:maxequalconvex}
Let $\{x_k\}$ be generated by \cref{algo:MIPP}. Then, for any $k \ge 1$, there exists $\theta_k^* := (\theta_{1,k}^*, \cdots, \theta_{m,k}^*)^\top \in \Delta^m$ such that the following three expressions hold:
\begin{equation}\label{eq:expression-1}
x_{k+1} \in \argmin_{z \in \mathbb{R}^n} \left\{ \sum_{i=1}^m \theta_{i,k}^* \big( f_i(z) - f_i(x_k) \big) + \frac{1}{2\lambda_k} \|z - y_k\|^2 \right\}.
\end{equation}
\begin{equation}
\max_{i\in[m]} \left( f_i(x_{k+1}) - f_i(x_k) \right) = \sum_{i=1}^m \theta_{i,k}^* \left( f_i(x_{k+1}) - f_i(x_k) \right).
\end{equation}
\begin{equation}\label{eq:expression-2}
\beta_k \sum_{i=1}^m \theta_{i,k}^* \partial f_i(x_{k+1}) + (x_{k+1} - 2x_k + x_{k-1}) + \frac{\alpha - 1}{k} (x_{k+1} - x_k) + \frac{1}{k} (x_k - x_{k-1}) \ni 0.
\end{equation}
\end{lemma}
\begin{proof}
According to \cref{eq:MIPP}, $x_{k+1}$ is the solution to the following subproblem:  
\begin{equation}\label{eq:expression-3}
\min_{z\in \mathbb{R}^n} \left\{G(z) = \Phi_{x_k}(z)+\frac{1}{2\lambda _{k}}\|z-y_k\|^2\right\}.
\end{equation}
This problem has an equivalent form:  
\begin{equation}\label{eq:equal-form-primal}
\begin{aligned}
	\min_{(z,\tau)\in \mathbb{R}^n\times \mathbb{R}} &\quad  \tau +\frac{1}{2\lambda _k}\|z-y_k\|^2,\\
	{\rm s.t.} &\quad f_i(z) -\tau \le f_i(x_k),\ i\in [m].
\end{aligned}
\end{equation}
For the Lagrangian $L(z,\tau,\theta) = \tau +\frac{1}{2\lambda _k}\|z-y_k\|^2+\sum_{ i = 1}^m \theta_{i}\big(f_i(z)-f_i(x_k)-\tau\big)$, consider the dual problem of \cref{eq:equal-form-primal}:  
\begin{equation}\label{eq:equal-form-dual}
\begin{aligned}
	\max_{\theta \in \mathbb{R}^m} &\quad \left\{\omega(\theta):=\min_{z\in \mathbb{R}^n,\tau \in \mathbb{R}}L(z,\tau,\theta)\right\}, \\
	{\rm s.t.}  &\quad \sum_{i=1}^m\theta_i  = 1,\quad \theta_i \ge 0, \ i\in[m]. 
\end{aligned}
\end{equation}
Let $\theta_k^*:=(\theta_{1,k}^*,\cdots,\theta_{m,k}^*)^\top  \in \Delta ^m$ be the solution to \cref{eq:equal-form-dual}. By the strong duality theorem \cite[Theorem A.1 and A.2]{beck2017first}, we obtain \cref{eq:expression-1} holds.  

Furthermore, according to \cref{eq:MIPP} and \cref{eq:expression-1}, we have  
\begin{equation*} 
\begin{aligned}&\sum_{i=1}^m \theta_{i,k}^*\big(f_i(x_{k+1})-f_i(x_k)\big)+\frac{1}{2\lambda_k}\|x_{k+1}-y_k\|^2 = \Phi_{x_k}(x_{k+1})+\frac{1}{2\lambda _k}\|x_{k+1}-y_k\|^2.
\end{aligned}
\end{equation*}
This proves \cref{eq:expression-2}.  

For \cref{eq:expression-1}, using the first-order optimality condition and noting that each $f_i$ has domain $\mathbb{R}^n$, we obtain  
\begin{equation*} 
\sum_{i = 1 }^m\theta _{i,k}^*\partial f_i(x_{k+1}) +\frac{1}{\lambda_k} (x_{k+1}-y_k)\ni0. 
\end{equation*}
Using the definitions of $y_k$ and $\lambda _k$, a straightforward calculation leads to \cref{eq:expression-3}.
\end{proof}

\subsection{Convergence rate}
In this section, define  
\begin{equation}\label{eq:sigma_k}
\sigma_k(z):=\min_{i\in [m]}(f_i(x_k)-f_i(z)).  
\end{equation}
and the auxiliary sequence  
\begin{equation}\label{eq:W_i}   
\mathcal W_{i,k}:= f_i(x_k)+\frac{1}{2\lambda_{k-1}}\|x_k-x_{k-1}\|^2.  
\end{equation}
\begin{lemma}  
Let $\{x_k\}$ be the sequence generated by \cref{algo:MIPP}. We have  
\begin{enumerate}[label=(\roman*)]
	\item $\{\mathcal W_{i,k}\}$ is non-increasing for $i\in[m]$.  
	\item $\{x_k\}\subseteq \mathcal L(F,F(x_0))$.
	\item  If each $f_i$ is bounded below, then $\lim_{k\to \infty }\mathcal W_{i,k}=\mathcal W_i^*\in \mathbb{R}$.  
\end{enumerate} 
\end{lemma}
\begin{proof} 
By \cref{lem:maxequalconvex}, there exists $\xi_k \in \partial \left(\sum_{i=1}^m\theta_{i,k}^* f_i(x_{k+1})\right)$ such that  
\begin{equation*} 
\xi_k +\frac{1}{\lambda _k}(x_{k+1}-y_k)=0,  
\end{equation*}
with $\xi_k =\sum_{i=1}^m \theta_{i,k}^*\xi_{i,k}$ and $\xi_{i,k}\in \partial f_i(x_{k+1})$. Then, using \cref{eq:beta_k}, we obtain  
\begin{equation}\label{eq:equal-inequal-1}   
\begin{aligned}  
\max_{i\in[m]}(f_i(x_{k+1})-f_i(x_k))&=\sum_{i=1}^m\theta_{i,k}^*(f_i(x_{k+1})-f_i(x_k))\le \frac1{\lambda_k}\langle y_k-x_{k+1},x_{k+1}-x_k\rangle\\  
&=\frac{1}{\lambda _k}\langle y_k-x_k,x_{k+1}-x_k\rangle  -\frac{1}{\lambda _k}\|x_{k+1}-x_k\|^2.  
\end{aligned}  
\end{equation}
Using the identity $\frac12\|a-b\|^2=\frac12\|a\|^2-\langle a,b\rangle +\frac12\|b\|^2$ together with $y_k =x_k+\frac{k-1}{k+\alpha -1}(x_k-x_{k-1})$, a straightforward computation gives  
\begin{equation} \label{eq:equal-inequal-2}   
\begin{aligned}  
\frac{1}{\lambda_k}\langle y_k-x_k,x_{k+1}-x_k\rangle &= \frac1{2\lambda _k}\|x_{k+1}-x_k\|^2+\frac{1}{2\lambda _{k-1}}\|x_k-x_{k-1}\|^2\\  
&\quad +\frac{1}{2\lambda_{k-1}}\left(\frac{\lambda _{k-1}}{\lambda_k}\frac{(k-1)^2}{(k+\alpha -1)^2}-1\right)\|x_k-x_{k-1}\|^2\\  
&\quad -\frac{1}{2\lambda_k}\|x_{k+1}-y_k\|^2.  
\end{aligned}  
\end{equation}   
From the definition of $\lambda_k$, we have  
\begin{equation}\label{eq:equal-inequal-3}   
\begin{aligned}  
\frac{\lambda_{k-1}}{\lambda _k}\frac{(k-1)^2}{(k+\alpha -1)^2}-1&=\frac{k(k-1)^2}{(k+\alpha-1)(k+\alpha -2)}-1<0.  
\end{aligned}  
\end{equation}  
Combining \cref{eq:equal-inequal-1,eq:equal-inequal-2,eq:equal-inequal-3},  yields  
\begin{equation*}   
\max_{i\in[m]}(f_i(x_{k+1})-f_i(x_k))+\frac1{2\lambda_k}\|x_{k+1}-y_k\|^2\le -\frac{1}{2\lambda_k}\|x_{k+1}-x_k\|^2+\frac{1}{2\lambda_{k-1}}\|x_k-x_{k-1}\|^2.  
\end{equation*}  
Since $f_i(x_{k+1})-f_i(x_k)\le \max_{i\in[m]}(f_i(x_{k+1})-f_i(x_k))$, the first claim follows.

(ii) From (i), we have  
\begin{equation*} 
f_i(x_k)+\frac1{2\lambda_{k-1}}\|x_k-x_{k-1}\|^2=\mathcal W_{i,k}\le\mathcal{W}_{i,1}=f_i(x_1)=f_i(x_0).  
\end{equation*}   
Thus the statement follows.  

(iii) The claim is immediate from (i).
\end{proof}
We define the following auxiliary sequence, termed the discrete Lyapunov-like function:
\begin{equation}\label{eq:desceret-Lyapunov-like} 
\mathcal E_k(z):=\underbrace{k^2\beta_k\sigma_k(z)}_{\textbf{Potential energy }\mathcal E_{k}^{\rm pot}(z) }+\underbrace{\frac{1}{2}\left\|(\alpha -1)(x_k-z)+(k-1)(x_{k}-x_{k-1})\right\|^2}_{\textbf{Mixed energy }\mathcal E_{k}^{\rm mix}(z) }.  
\end{equation}  
For the mixed energy, we first perform the following computation:  
\begin{equation}\label{eq:inequal-mixed-1}   
\begin{aligned}  
	&[(\alpha -1)(x_{k+1}-z) +k(x_{k+1}-x_k)] -[(\alpha -1)(x_k-z)+(k-1)(x_k-x_{k-1})]\\  
	&=(\alpha -1) (x_{k+1}-x_k)+(x_k-x_{k-1})+k(x_{k+1}-2x_k+x_{k-1})\\  
	&=-k\beta_k \left(\frac{1}{\lambda_{k}}(x_{k+1}-y_k)\right)\\  
	&=-k\beta_k \xi_k.  
\end{aligned}  
\end{equation}
where $\xi _k =\sum_{i = 1}^m \theta_{i,k}^* \xi_{i,k}$ and $\xi_{i,k}\in \partial f_i(x_{k+1})$. Using this, and the identity $\frac12\|a\|^2-\frac12\|b\|^2 +\frac{1}{2}\|a-b\|^2= \langle a-b,a\rangle$, we have  
\begin{equation}\label{eq:inequal-mixed-2}
\begin{aligned}  
	&\mathcal E_{k+1}^{\rm mix}(z)-\mathcal E_k^{\rm mix}(z)+\frac12\|k\beta_k \xi_k\|^2\\  
	&   =\left\langle -k\beta_k \xi_k,(\alpha -1)(x_{k+1}-z)+k(x_{k+1}-x_{k})\right\rangle\\  
	&\le -~(\alpha -1) k \beta_k\sum_{i=1}^m  \theta_{i,k}^*\left(f_i(x_{k+1})-f_i(z)\right) - k^2 \beta_k\sum_{i=1}^m  \theta_{i,k}^*(f_i(x_{k+1})-f_i(x_k))\\  
	&\le -~(\alpha -1) k\beta_k \sigma_{k+1}(z)-k^2\beta_k\max_{i\in[m]} \left(f_i(x_{k+1})-f_i(x_k)\right).  
\end{aligned}  
\end{equation}  
where the last inequality follows from $\min_{i\in[m]}(f_i(x_{k+1})-f_i(z))\le \sum_{i = 1}^m \theta_{i,k}^*  (f_i(x_{k+1})-f_i(z))$ and \cref{lem:maxequalconvex}. 
For the potential energy, we have  
\begin{equation}\label{eq:inequal-potential}
\begin{aligned}  
	&\mathcal E_{k+1}^{\rm pot}(z)-\mathcal E_k^{\rm pot}(z)\\  
	& =  \left[(k+1)^2(\beta_{k+1}-\beta_k)+(2k+1)\beta_k\right]\sigma_{k+1}(z)+k^2\beta_k(\sigma_{k+1}(z)-\sigma_k(z))\\  
	&\le \left[(k+1)^2(\beta_{k+1}-\beta_k)+(2k+1)\beta_k\right]\sigma_{k+1}(z)+k^2\beta_k\Phi_{x_k}(x_{k+1}).
\end{aligned}  
\end{equation}  
where $\Phi_{x_k}(x_{k+1}) = \max_{i\in[m]}(f_i(x_{k+1})-f_i(x-k))$. Based on \cref{eq:inequal-mixed-2,eq:inequal-potential}, we can prove the following theorem.  

\begin{theorem}\label{thm:convergence-rate} 
Suppose that $f_i$ is lower semicontinuous for all $i\in[m]$.  Let $\{x_k\}$ be the sequence generated by \cref{algo:MIPP} with initial point $x_0$. Then,
\begin{enumerate}[label=(\roman*)]
	\item $u_0(x_k ) = O(1/k^2\beta_k)$;
	\item $\|x_k-x_{k-1}\|=O(1/k)$;
\end{enumerate}  
\end{theorem}
\begin{proof}
From \cref{eq:inequal-mixed-2,eq:inequal-potential}, we obtain  
\begin{equation*} 
 \begin{aligned}  
	&\mathcal E_{k+1}(z)-\mathcal E_k(z)+\frac{1}{2}\|k\beta_k\xi_k\|^2\\  
	& \le   \left[(k+1)^2(\beta_{k+1}-\beta_k)+(2k+1)\beta_k\right]\sigma_{k+1}(z)-(\alpha -1) k\beta_k \sigma_{k+1}(z)\\  
	&=\left[(k+1)^2\beta_{k+1}-k(k+\alpha -1)\beta_k\right]\sigma_{k+1}(z).  
\end{aligned}  
\end{equation*}   
From \cref{eq:beta_k}, for any $l\ge 1$, we get  
\begin{equation*}   
 \mathcal E_{l+1}(z)-\mathcal E_{l}(z)+\frac12\|l\beta_l \xi_l\|^2\le 0.  
\end{equation*}
Summing from $l=1$ to $l= k -1$ gives  
\begin{equation}\label{eq:inequal-lyapunov-bounded}   
\mathcal E_{k}(z) +\frac12\sum_{l=1}^{k-1}\|l\beta_l\xi_l\|^2\le \mathcal E_1(z)\le \beta_1 u_0(x_0)+\frac{(\alpha -1)^2}2\|x_0-z\|^2.  
\end{equation}   
From \cref{eq:inequal-lyapunov-bounded} we can separately derive the two statements of the theorem. Specifically,  
\begin{itemize}
\item 	Since $\mathcal E_k(z)\ge \sigma_k(z)$, and by \cref{thm:sup-inf}, we have  
\begin{equation*}
\begin{aligned}   
 k^2\beta_ku_0(x_k)&\le k^2\beta_k\sup_{F^* \in F\left(\mathcal L\mathcal{P}_w (F, F(x_0))\right)} \inf_{x \in F^{-1}(F^*)} \sigma_k(z)+\frac{1}{2}\sum_{l=1}^{k-1}\|\ell \beta_{l}\xi_{l}\|^2\\&\le \beta_1u_0(x_0)+(\alpha -1)^2\mathbf{R}.  
\end{aligned}
\end{equation*}   
This proves (i).  
\item Taking $z =x_k$, we obtain  
\begin{equation*}   
\frac12(k-1)^2\|x_k-x_{k-1}\|^2\le \beta_1u_0(x_0)+(\alpha -1)^2\left(\|x_0\|^2+\sup_{k\ge 1}\|x_k\|^2\right).  
\end{equation*}  
This gives part (ii). 
\end{itemize}
The proof is complete.
\end{proof}
\subsection{Convergence analysis}
The following provides the convergence analysis of the iterative sequence.
\begin{lemma}
	Suppose that $f_i$ is lower semicontinuous for all $i\in[m]$, and \cref{assume:boundedsup} hold true. Let  $\{x_k\}$ be a bounded sequence generated by \cref{algo:MIPP} with initial point $x_0$. Assuming $f_i$ is bounded from below for $i\in[m]$. Then, 
	\begin{enumerate}[label=(\roman*)]
		\item $\lim_{k\to \infty }f_i(x_k)=f_i^\infty $ exists;
		\item For any cluster point $z^*$ of $\{x_k\}$, the limit  
		  \begin{equation}    
		  \lim_{k\to \infty }\mathcal E_k(z^*) = \mathcal E_{z^*}^\infty \in \mathbb{R},
		  \end{equation}  
		   exists.  
	\end{enumerate} 
\end{lemma}
\begin{proof}
(i) By the non-increase of $\mathcal W_{i,k}$, for any $k\ge l\ge 1$,  
\begin{equation*} 
\mathcal W_{k}:=u_0(x_k)+\frac{1}{2\lambda_{k-1}}\|x_k-x_{k-1}\|^2\le u_0(x_l)+\frac{1}{2\lambda_{l-1}}\|x_{l}-x_{l-1}\|^2.  
\end{equation*}   
Hence $\{\mathcal W_k\}$ is non-increasing. Since $u_0(x_k)\ge 0$, the limit $\lim_{k\to \infty }\mathcal W_k=:\mathcal W^*\in \mathbb{R}_+$ exists. Moreover, because $u_0(x_k)\to 0$ as $k\to \infty$, we have $\frac{1}{2\lambda_{k-1}}\|x_k-x_{k-1}\|^2\to \mathcal W^*$ as $k\to \infty$. Consequently,  
\begin{equation*}   
f_i^\infty :=\lim_{k\to \infty }f_i(x_k)=\lim_{k\to \infty}\left(\mathcal W_{i,k}-\frac{1}{2\lambda_{k-1}}\|x_k-x_{k-1}\|^2\right)\in \mathbb{R}  
\end{equation*}   
exists.  

(ii) Since each $f_i$ is lower semicontinuous and $z^*$ is a cluster point of $\{x_k\}$,  
\begin{equation*}   
\begin{aligned}  
	\sigma_k(z^*) &= \min_{i\in[m]}\left(\mathcal W_{i,k}-\frac{1}{2\lambda_{k-1}}\|x_k-x_{k-1}\|^2-f_i(z^*)\right)\\
	&\ge \min_{i=1,\cdots,m}\left(f_i^\infty -f_i(z^*)\right)-\frac{1}{2\lambda_{k-1}}\|x_k-x_{k-1}\|^2\\
	&\ge -\frac1{2\lambda_{k-1}}\|x_k-x_{k-1}\|^2\\
	&=-\frac{1}{2}\frac{k+\alpha -2}{(k-1)\beta_{k-1}}\|x_k-x_{k-1}\|^2,
\end{aligned}  
\end{equation*}   
where the first inequality uses $\mathcal W_{i.k}\ge \lim_{k\to \infty }\mathcal W_{i,k}\ge f_i^\infty$, the second follows from the lower semicontinuity of $f_i$, and the last equality comes from the definition of $\lambda _k$. Therefore,  
\begin{equation*}   
k^2\beta_k\sigma_k(z^*) \ge -\frac12\frac{k+\alpha -2}{k-1} \frac{\beta_k}{\beta_{k-1}} k^2 \|x_k-x_{k-1}\|^2=-\frac{1}{2}({k+\alpha -2})^2\|x_k-x_{k-1}\|^2. 
\end{equation*}   
The last equality follows from \cref{eq:beta_k}. By \cref{thm:convergence-rate}, the right-hand side is bounded below by a finite constant, so there exists $\widetilde {\mathbf{R}}\in \mathbb{R}$ such that  
\begin{equation*}   
 k^2\beta_k \sigma_k(z^*)\ge\widetilde {\mathbf{R}}>-\infty.  
\end{equation*}   
Thus $\mathcal E_{k}(z^*)$ is non-increasing and bounded below, hence its limit exists.  
\end{proof} 
{The following lemma is a discrete version of \cref{lem:continuous-limit}.}
\begin{lemma}[{\cite[Lemma A.5]{boct2025fast}}]\label{lem:q_k}
    Let $ a \geq 1 $ and $\{q_k\}_{k \geq 0}$ be a bounded sequence in $\mathbb{R}$ such that
    \begin{equation*} 
    \lim_{k \to +\infty} \left( q_{k+1} + \frac{k}{a} (q_{k+1} - q_k) \right) = l \in \mathbb{R}.
    \end{equation*} 
    Then, it holds $\lim_{k \to +\infty} q_k = l$.
\end{lemma}
Under certain conditions, we prove the point convergence of \cref{algo:MIPP} as follows:
\begin{theorem}
Suppose that $f_i$ is lower semicontinuous for all $i\in[m]$, and \cref{assume:boundedsup} hold true. Let $\{x_k\}$ be the bounded sequence generated by \cref{algo:MIPP}. If one of the following two conditions holds: 
\begin{enumerate}[label=(\alph*)]
	\item  Each $f_i$ is continuous for $i\in [m]$;
	\item  $F^\infty =(f_1^\infty ,\cdots,f_m^\infty )^\top \in F(\mathcal P)$.  
\end{enumerate} 
Then $x_k \to x^*$ as $k\to \infty$.  
\end{theorem} 
\begin{proof}
 By the lower semicontinuity  of $f_i$, for any cluster point $z^*$ of $\{x_k\}$,  
\begin{equation*}   
f_i^\infty =\lim_{k\to \infty }f_i(x_k)\ge \liminf_{x\to z^*} f_i(x)\ge  f_i(z^*).  
\end{equation*}   
Because $F^\infty \in F(\mathcal P)$, there is no point $y$ satisfying $f_i(y)\le f_i^\infty$ for all $i\in[m]$ and $f_j(y)<f_j^\infty$ for some $j$. Hence  
\begin{equation*}   
f_i^\infty = f_i(z^*).  
\end{equation*} 
For any two cluster points $z_1$ and $z_2$ of $\{x_k\}$,  
\begin{equation*}   
\begin{aligned}  
	\mathcal E_k(z_j)&= k^2\beta_k \sigma_k(z_j)+\frac12\left\|(\alpha -1)(x_k-z_j)+(k-1)(x_k-x_{k-1})\right\|^2\\
	&=k^2 \beta_k \sigma_{\infty } +\frac{(k-1)^2}{2}\|x_k-x_{k-1}\|^2
	+ \frac{(\alpha -1)^2}{2}\|x_k-z_j\|^2\\
	&\qquad +\ (\alpha -1)(k-1)\langle x_k-z_j,x_k-x_{k-1}\rangle,
\end{aligned}  
\end{equation*}   
and  
\begin{equation*}   
\begin{aligned}  
	\mathcal E_{k}(z_1)-\mathcal E_k(z_2)&=(\alpha -1)^2\left(\frac12\|x_k-z_1\|^2-\frac12\|x_k-z_2\|^2\right)\\
	&\qquad +\ (\alpha -1)(k-1)\langle z_2-z_1,x_k-x_{k-1}\rangle .
\end{aligned}  
\end{equation*}   
Notice that  
\begin{equation*}   
\begin{aligned}  
	\|x_k-z_1\|^2-\|x_k-z_2\|^2 &=\|x_k-z_2\|^2+2\langle x_k-z_2,z_2-z_1\rangle +\|z_2-z_1\|^2-\|x_k-z_2\|^2\\
	&=2\langle x_k-z_2,z_2-z_1\rangle +\|z_2-z_1\|^2,\\
	\|x_{k-1}-z_1\|^2-\|x_{k-1}-z_2\|^2 &= 2\langle x_{k-1}-z_2,z_2-z_1\rangle +\|z_2-z_1\|^2.
\end{aligned}  
\end{equation*}   
If we set  
\begin{equation*}   
h_k :=\frac12\|x_k-z_1\|^2-\frac12\|x_k-z_2\|^2,  
\end{equation*}  
then  
\begin{equation*}  
\begin{aligned}  
	\mathcal E_{k}(z_1)-\mathcal E_k(z_2)&=(\alpha -1)^2h_k+2 (\alpha -1)(k-1)\cdot (h_k-h_{k-1}) .
\end{aligned}  
\end{equation*}   
Hence  
\begin{equation*}   
\frac{1}{(\alpha -1)^2}(\mathcal E_k(z_1)-\mathcal E_k(z_2)) = h_k +\frac{k-1}{\frac{1}{2}(\alpha -1)}(h_k-h_{k-1}).  
\end{equation*}  
Since the left-hand side has a limit, the right-hand side also has a limit, by \cref{lem:q_k}, we have
\begin{equation*}   
\lim_{k\to \infty } h_k =\ell.  
\end{equation*}   
Suppose $\{\bar x_k\}$ and $\{\tilde x_k\}$ are subsequences of $\{x_k\}$ with $\bar x_k\to z_1$ and $\tilde x_k \to z_2$ as $k\to \infty$. Then  
\begin{equation*}   
\begin{aligned}  
	2\ell = \lim_{k\to \infty}\left( \|\bar x_k -z_1\|^2-\|\bar x_k -z_2\|^2\right)&=-\|z_1-z_2\|^2,\\
	2\ell =\lim_{k\to \infty}\left(\|\tilde  x_k-z_1\|^2-\|\tilde x_k -z_2\|^2\right)&=\|z_1-z_2\|^2.
\end{aligned}  
\end{equation*}   
Therefore $z_1=z_2=:x^*$. That is, $\lim_{k\to \infty }x_k =x^*$. If (a) holds, $x^*$ is a weakly Pareto optimal solution; if (b) holds, $x^*$ is a Pareto optimal solution.
\end{proof}  
\section{Numerical experiments}\label{sec:num}
This section presents the performance of the \cref{eq:MIPP} across different problems through numerical experiments. All numerical experiments were performed in the MATLAB 2021a environment on a personal computer with an Intel(R) Core(TM) i5-8300H CPU @ 2.30GHz   2.30 GHz processor and 8GB of RAM. We set the termination condition tolerance to $\varepsilon =1e-4$ and the maximum number of iterations to $k_{\max} =1000$.  
\begin{table}[h]
    \centering
    \caption{Smooth test problems.}
    \label{tab:test_problems-smooh}
    \begin{tabular}{c c c c c c}
        \toprule
        \textbf{Problem} & $\boldsymbol{n}$ & $\boldsymbol{m}$ & $\boldsymbol{x_L}$ & $\boldsymbol{x_U}$ & \textbf{Reference} \\
        \midrule
        JOS1        & 20 & 2 & $(-15,\dots,-15)$ & $(15,\dots,15)$      & \cite{mita2019nonmonotone} \\
        SD        & 4 & 2 & $(1,\sqrt 2,\sqrt 2.1)$       & $(3,3,3,3)$              & \cite{mita2019nonmonotone} \\
        TOI4       & 4 & 2 & $(-2,-2,-2,-2)$         & $(5,5,5,5)$              & \cite{mita2019nonmonotone} \\
        FDS        & 10 & 3 & $(-2,\dots,-2)$   & $(2,\dots,2)$        & \cite{mita2019nonmonotone} \\
        TRIDIA        & 3 & 3 & $(-1,-1,-1)$         & $(1,1,1)$              & \cite{mita2019nonmonotone} \\
        SP1       & 20 & 3 & $(-15,\cdots,-15)$           & $(15,\cdots,15)$              & \cite{sonntag2024fast} \\
        \bottomrule
    \end{tabular}
\end{table}
\begin{table}[h]
	\centering
	\caption{Non-smooth test problems.}
	\label{tab:test_problems-nonsmooth}
	\begin{tabular}{c c c c c c}
		\toprule
		\textbf{Problem} & $\boldsymbol{n}$ & $\boldsymbol{m}$ & $\boldsymbol{x_L}$ & $\boldsymbol{x_U}$ & \textbf{Reference} \\
		\midrule
		\multirow{2}{*}{MKM1} & 2 & 2 & $(-5,-5)$         & $(5,5)$              & \cite{montonen2018multiple} \\
		&\multicolumn{4}{c}{$F(x) = \begin{pmatrix}
				\max\left\{x_1^2 + (x_2 - 1)^2,  (x_1 + 1)^2\right\} \\
				\max\left\{2x_1 + 2x_2,  x_1^4 + x_2^2\right\}
			\end{pmatrix}$}\\ 
		\midrule
		\multirow{2}{*}{MKM2} & 2 & 2 & $(-5,-5)$         & $(5,5)$              & \cite{montonen2018multiple} \\
			&\multicolumn{4}{c}{$F(x) = \begin{pmatrix}
					\max\left\{(x_1 - 2)^2 + (x_2 + 2)^2,  x_1^2 + 8x_2\right\} \\
					\max\left\{5x_1 + x_2,  x_1^2 + x_2^2\right\}
				\end{pmatrix} $}\\
		\midrule 
		\multirow{2}{*}{MKM3}        & 2 & 2 & $(-5,-5)$         & $(5,5)$              & \cite{montonen2018multiple}  \\
		&\multicolumn{4}{c}{$	F(x) = \begin{pmatrix}
				\max\left\{x_1^4 + x_2^2,  (2 - x_1)^2 + (2 - x_2)^2,  2e^{x_2 - x_1}\right\} \\
				\max\left\{-x_1 - x_2,  -x_1 - x_2 + x_1^2 + x_2^2 - 1\right\}
			\end{pmatrix} $}\\
		\bottomrule
	\end{tabular}
\end{table}

We randomly selected 200 points for each test problem within the corresponding domain according to \cref{tab:test_problems-smooh} and \cref{tab:test_problems-nonsmooth}. For each choice of $\alpha \in \{3, 5, 10, 25\}$, the iterative sequence was generated using \cref{algo:MIPP}. The contents and main results of the numerical experiments are as follows:

\begin{itemize}
    \item \cref{fig:smooth-pareto-front,fig:nonsmooth-pareto-front} display the weak Pareto fronts of smooth and nonsmooth bi-objective test problems, respectively. It can be observed that for smooth functions, different choices of $\alpha$ have little influence on the weak Pareto front produced by the algorithm; for nonsmooth functions, a larger $\alpha$ leads to a more uneven weak Pareto front generated by the algorithm, a phenomenon that is more pronounced on test problem MKM3 (\cref{fig:MKM3}). 
    \item \cref{tab:convexresult} records the total number of iterations and the iteration time for 200 points under different choices of $\alpha$. Each algorithm was evaluated using Performance Evaluation \cite{dolan2002benchmarkingperformance}, as shown in \cref{fig:performance}. The experimental results indicate that as $\alpha$ increases, both the convergence time and the iteration time of the point sequence generated by \cref{algo:MIPP} decrease. This outcome supports the correctness of \cref{thm:convergence-rate} .
\end{itemize} 

\begin{figure}
	\centering
	\begin{subfigure}[b]{0.32\textwidth}
		\centering
		\includegraphics[width=\textwidth]{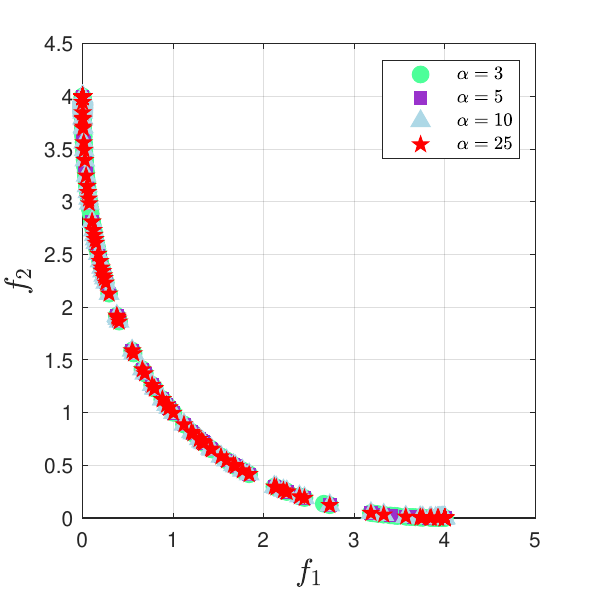}
		\caption{JOS1}
		\label{fig:JOS1}
	\end{subfigure}
	\hfill     
	\begin{subfigure}[b]{0.32\textwidth}
		\centering
		\includegraphics[width=\textwidth]{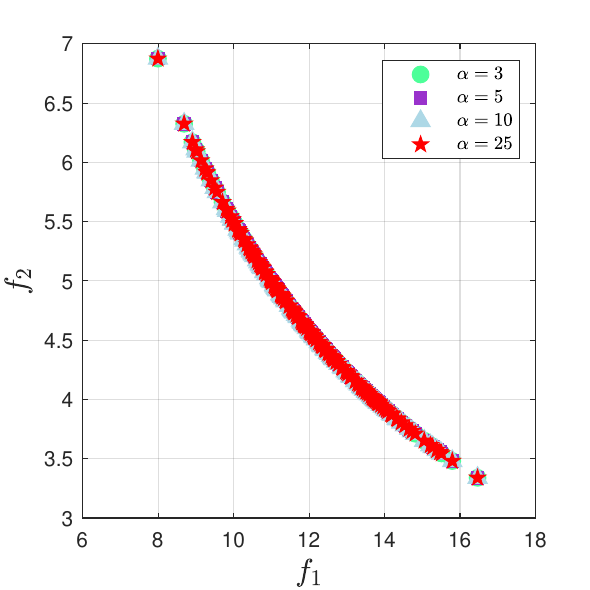}
		\caption{SD}
		\label{fig:SD}
	\end{subfigure}
	\hfill     
	\begin{subfigure}[b]{0.32\textwidth}
		\centering
		\includegraphics[width=\textwidth]{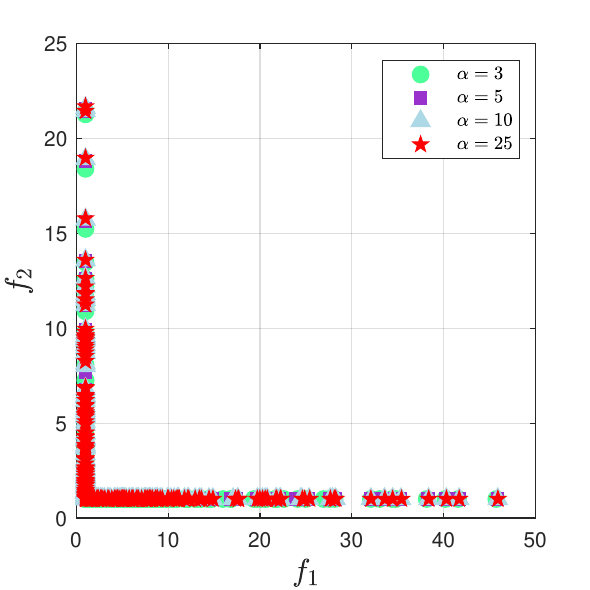}
		\caption{TOI4}
		\label{fig:TOI4}
	\end{subfigure}
	\caption{Smooth test problems.}
	\label{fig:smooth-pareto-front}
\end{figure}
\begin{figure}
	\centering
	\begin{subfigure}[b]{0.32\textwidth}
		\centering
		\includegraphics[width=\textwidth]{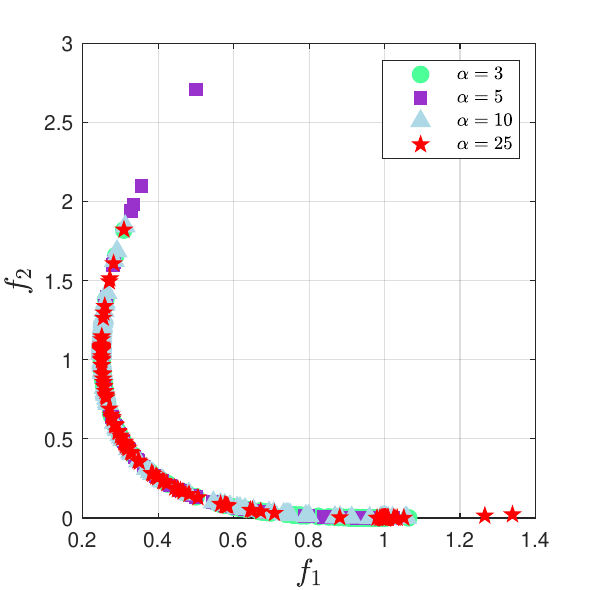}
		\caption{MKM1}
		\label{fig:MKM1}
	\end{subfigure}
	\hfill     
	\begin{subfigure}[b]{0.32\textwidth}
		\centering
		\includegraphics[width=\textwidth]{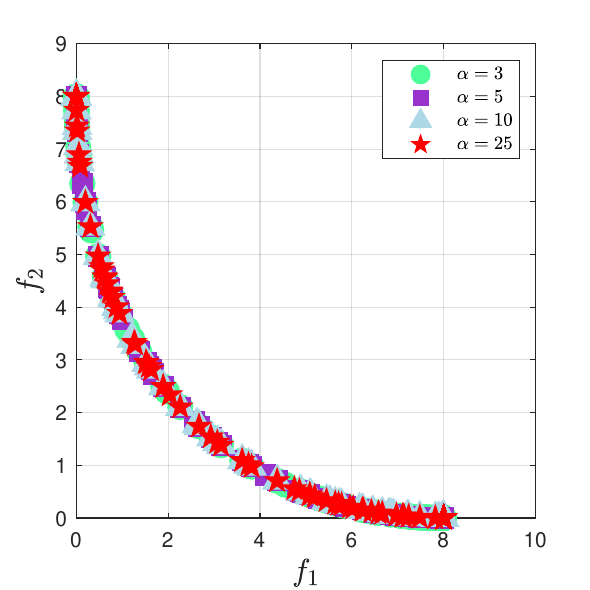}
		\caption{MKM2}
		\label{fig:MKM2}
	\end{subfigure}
	\hfill     
	\begin{subfigure}[b]{0.32\textwidth}
		\centering
		\includegraphics[width=\textwidth]{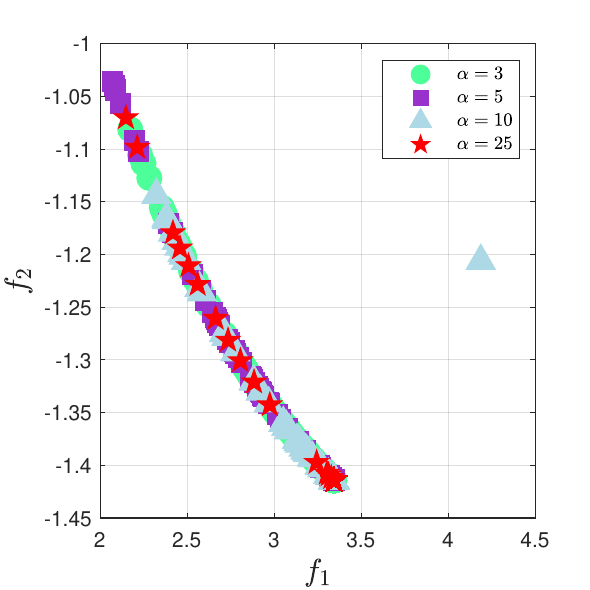}
		\caption{MKM3}
		\label{fig:MKM3}
	\end{subfigure}
	\caption{Non-smooth test problems.}
	\label{fig:nonsmooth-pareto-front}
\end{figure}
\begin{table}[h]
    \centering
    \caption{Total iterations (Iter) and Total CPU time (Time (s)) of tested algorithms implemented on different test problems}
    \label{tab:convexresult}
    \begin{tabular}{cccccccccc}
        \toprule
        \multirow{2}{*}{Problem} & \multicolumn{2}{c|}{$\alpha =3$} & \multicolumn{2}{c|}{$\alpha=5$} & \multicolumn{2}{c|}{$\alpha=10$} & \multicolumn{2}{c|}{$\alpha=25$} \\
        \cmidrule(lr){2-3}\cmidrule(lr){4-5}\cmidrule(lr){6-7}\cmidrule(lr){8-9}
        & {\textbf{Time}} & \textbf{Iter} & \textbf{Time} & \textbf{Iter} & \textbf{Time} & \textbf{Iter} & \textbf{Time} & \textbf{Iter} \\
        \midrule
        JOS1 & 506.73 & 21844 & 115.22 & 4243 & 60.50 & 2200 & 38.91 & 1400 \\
        SD & 118.27 & 7073 & 41.45 & 2557 & 27.53 & 1553 & 20.00 & 1098 \\
        TOI4 & 64.64 & 3900 & 33.11 & 2005 & 22.38 & 1357 & 18.47 & 1080 \\
        FDS & 376.23 & 16551 & 98.22 & 3862 & 54.44 & 1872 & 41.84 & 1259 \\
        TRIDIA & 47.83 & 2224 & 28.69 & 1397 & 21.39 & 1092 & 19.78 & 899 \\
        SP1 & 1113.92 & 37650 & 190.98 & 5370 & 107.56 & 2282 & 77.39 & 1488 \\
        MKM1 & 81.55 & 2140 & 47.66 & 1791 & 35.58 & 1320 & 30.08 & 1101 \\
        MKM2 & 41.19 & 1684 & 33.31 & 1470 & 27.89 & 1227 & 24.39 & 1098 \\
        MKM3 & 65.44 & 3564 & 42.98 & 1971 & 35.36 & 1357 & 26.20 & 1120 \\
        \bottomrule
    \end{tabular}
\end{table}

\begin{figure}[H]
	\centering
	\begin{subfigure}[b]{0.48\textwidth}
		\centering
		\includegraphics[width=\textwidth]{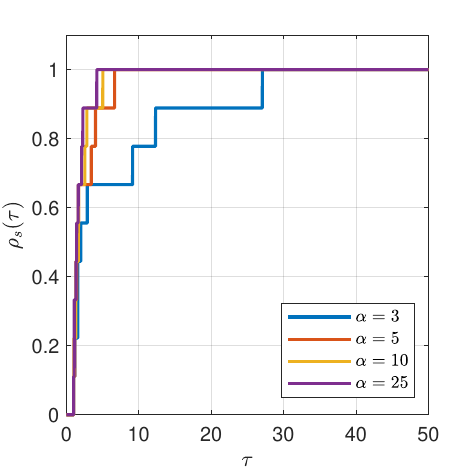}
		\caption{Time}
		\label{fig:subi1}
	\end{subfigure}
	\hfill     
	\begin{subfigure}[b]{0.48\textwidth}
		\centering
		\includegraphics[width=\textwidth]{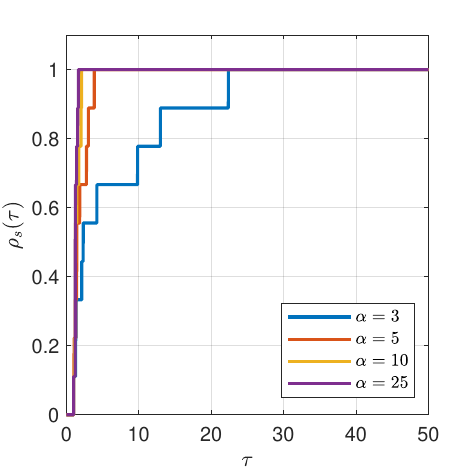}
		\caption{Iteration}
		\label{fig:subi2}
	\end{subfigure}
	\caption{Performance evaluation}
	\label{fig:performance}
\end{figure}
\section{Conclusion}\label{sec:con}
In this work, we introduced and analyzed a multiobjective inertial gradient system with time scaling \cref{eq:MITS} and its discretized counterpart, the multiobjective inertial proximal point method \cref{eq:MIPP}. We established the existence of solution trajectories for \ref{eq:MITS} and proved, using Lyapunov analysis, that the trajectory achieves an arbitrarily fast sublinear convergence rate of $O(1/t^{2}\beta(t))$ with respect to a merit function, along with convergence to a weakly Pareto optimal solution. The corresponding discrete algorithm \ref{eq:MIPP} inherits the fast $O(1/k^{2}\beta_{k})$ rate and converges under suitable conditions. Numerical experiments validate the theoretical findings and demonstrate the efficiency of the proposed method. Our results extend and accelerate existing inertial dynamics and proximal schemes to the multiobjective setting.
\begin{appendices}
\appendix
\section{Poof of Existence Theorem}\label{appendix:existence}
In this appendix, we prove the existence of solutions for the general multiobjective gradient flow. Consider the following Cauchy problem:
\begin{equation}\label{eq:CP}
\left\{\begin{aligned} &\alpha (t)\dot x(t)+\proj_{B(t,x(t),\dot x(t))+\ddot x(t)}(0) = 0,~~\text{for all} \ t>t_0,\\
&x(t_0) = x_0, ~~~\dot x(t_0) =v_0, 
\end{aligned}\right.\tag{CP}
\end{equation}
where $\alpha(t)$ is a continuous function from $[t_0,+\infty )$ to $\mathbb{R}_{++}$, and $B$ is a set-valued mapping:
\begin{equation}\label{eq:set-B}
B:(t_{\min},+\infty)\times\mathbb{R}^n\times \mathbb{R}^n \rightrightarrows \mathbb{R}^n,
\end{equation}
where $t_{\min}<t_0$. We now give the definition of the solution of \cref{eq:CP}.

\begin{definition} \label{def:soultion-general}
	{ A function $x:[t_0,+\infty )\to \mathbb{R}^n$, $t\mapsto x(t) $, is called a solution  of \cref{eq:CP} if it satisfies the following conditions: } 
	\begin{enumerate}[label=(\roman*)] 
		\item $x(\cdot)\in C^1([t_0,+\infty ))$.  
		\item $\dot x(t)$ is absolutely continuous on $[t_0,T]$ for any $T\ge t_0$.  
		\item There exists a measurable function $\ddot x(t):[t_0,+\infty )\to \mathbb{R}^n$ such that $\dot x(t)=\dot x(t_0)+\int_{t_0}^t\ddot x(s)ds$ for $t\ge t_0$.  
		\item $\dot x(\cdot) $ is almost everywhere differentiable and $\frac{d}{dt}\dot x(t)=\ddot x(t)$ holds for almost all $t\in [t_0,+\infty )$.  
		\item $\alpha (t)\dot x(t)+\proj_{B(t,x(t),\dot x(t))+\ddot x(t)}(0) = 0$ for almost all $t\in[t_0,+\infty )$.
		\item $x(t_0)=x_0$ and $\dot x(t_0)=v_0$.  
	\end{enumerate}
\end{definition}
 We refer to the following assumption as the Attouch-Goudou conditions, since their core components were introduced by these two authors in \cite[Section 3]{Attouch2014}.
\begin{definition}
For two closed convex subsets $C$ and $D$ in $\mathbb{R}^n$, the Hausdorff distance is defined as  
\begin{equation*} 
d_H(C, D) = \max \left\{ \sup_{x \in C} d(x, D), \, \sup_{x \in D} d(x, C) \right\}.
\end{equation*} 
Equivalently (see \cite{dontchev2009implicit}),  
\begin{equation*} 
d_H(C, D) = \sup_{x \in \mathbb{R}^n} \bigl| d(x, C) - d(x, D) \bigr|.
\end{equation*} 
\end{definition}
\begin{definition}
    Let $\Omega \subseteq \mathbb{R}^{d}$ be an open set. We say that a set-valued mapping $Q:\Omega \rightrightarrows \mathbb{R}^{p}$ is upper semi-continuous on $\Omega$ if  
    \begin{enumerate}[label=(\roman*)] 
        \item $Q(x)$ is compact for any $x \in \Omega$; 
        \item for any $x_0 \in \Omega$ and any $\varepsilon > 0$, there exists $\delta > 0$ such that  
        $$Q(x) \subseteq Q(x_0) + \mathcal{B}_{\varepsilon}(0), \quad \text{for all } x \in \mathcal{B}_{\delta}(x_0).$$
    \end{enumerate} 
\end{definition}
\begin{assum}\label{assum:A-G-condition}
Let $B$ be the set-valued mapping defined in \cref{eq:set-B} and satisfy the following conditions:
\begin{enumerate}[label=(\roman*)] 
\item[$\rm (i)$] For every $(t,u,v) \in (t_{\min}, +\infty) \times \mathbb{R}^n \times \mathbb{R}^n$, the set $B(t,u,v)$ is compact and convex.
\item[$\rm (ii)$] $B$ is Hausdorff continuous; that is, for every $(\bar t,\bar u,\bar v)$,
\begin{equation*} 
\begin{aligned}
&~~~~\lim_{(t,u,v) \to (\bar t,\bar u,\bar v)} d_{H}\left(B(t,u,v), B(\bar t,\bar u,\bar v)\right)  = 0.
\end{aligned}
\end{equation*} 
\item[$\rm (iii)$] For every $(t,u,v) \in (t_{\min}, +\infty) \times \mathbb{R}^n \times \mathbb{R}^n$, the following inequality holds:
\begin{equation*} 
\sup_{x \in B(t,u,v)} \|x\| \le \varphi(t) \left(1 + \|(u,v)\|\right),
\end{equation*} 
where $\varphi : (t_{\min}, +\infty) \to \mathbb{R}_{+}$ is a function that is bounded on any closed interval.
\end{enumerate}
\end{assum}
In \cite{YingdongYin2026BalancedGradientFlow}, the third condition of the Attouch-Goudou conditions given by Yin et al. is:
\begin{equation*} 
\sup_{x \in B(t,u,v)} \|x\| \le c(1 + \|(t, u, v)\|).
\end{equation*} 
The condition specified in this paper is strictly weaker than the one above, as evidenced by the following inequalities:
\begin{equation*} 
\sup_{x \in B(t,u,v)} \|x\| \le c(1 + \|(t, u, v)\|) \le c(1 + |t| + \|(u, v)\|) \le c(1 + |t|)(1 + \|(u, v)\|).
\end{equation*} 
\end{appendices}
\begin{example}
Assume that $ f_i $, $ i = 1\in [m]$, have Lipschitz continuous gradients. For $ \delta > 0 $, define the set-valued mapping $ \widetilde{C}: (t,x) \mapsto \beta(t) C(x) $. Then $ \widetilde{C} $ satisfies \cref{assum:A-G-condition}.
\end{example}
\begin{proof}
(i) Since the convex hull of a finite number of vectors is a closed set, this conclusion is evidently valid.

(ii) For any fixed $(\bar t,\bar x) \in [t_0, +\infty) \times \mathbb{R}^n$, based on the definition of the Hausdorff distance, we have for any $(t,x) \in [t_0, +\infty) \times \mathbb{R}^n$:
\begin{equation}
d_{H}(\beta(t) C(x), \bar \beta(t) C(\bar x)) \le d_{H}(\beta(t) C(x), \bar \beta(t) C(x)) + d_{H}(\bar \beta(t) C(x), \bar \beta(t) C(\bar x)).
\end{equation}
Note that for any $\xi \in \mathbb{R}^n$:
\begin{equation*} 
\begin{aligned}
d(\xi, \beta(t) C(x)) &= \inf_{a \in C(x)} \|\xi - \beta(t) a\| = \inf_{a \in C(x)} \|\xi - \bar \beta(t) a + (\bar \beta(t) - \beta(t))a\| \\
&\le \inf_{a \in C(x)} \|\xi - \bar \beta(t) a\| + \sup_{a\in C(x)}\|a\| \cdot |\beta(t) - \bar \beta(t)| \\&= d(\xi, \bar \beta(t) C(x)) + \sup_{a\in C(x)}\|a\| \cdot |\beta(t) - \bar \beta(t)|.
\end{aligned}
\end{equation*} 
Therefore, the following inequality holds:
\begin{equation}
\begin{aligned}
d_{H}(\beta(t) C(x), \bar \beta(t) C(x)) &= \sup_{\xi \in \mathbb{R}^n} \left|d(\xi, \beta(t) C(x)) - d(\xi, \bar \beta(t) C(x))\right|\le \sup_{a \in C(x)} \|a\| \cdot |\beta(t) - \bar \beta(t)| \\
&= \left\| \sum_{k=1}^m \lambda_i \nabla f_i(x) \right\| \cdot |\beta(t) - \bar \beta(t)| \\
&\le \left(L \cdot \|x - \bar x\| + \max_{i\in[m]} \|\nabla f_i(\bar x)\| \right) \cdot |\beta(t) - \bar \beta(t)|.
\end{aligned}
\end{equation}
According to \cite[Lemma 3.2]{Attouch2014}, we obtain:
\begin{equation}
\begin{aligned}
d_{H}(\bar \beta(t) C(x), \bar \beta(t) C(\bar x)) &= \bar \beta(t) \cdot d_{H}(C(x), C(\bar x)) \le \bar \beta(t) L \|x - \bar x\|.
\end{aligned}
\end{equation}
Consequently, for any sequence $(t_n, x_n) \to (\bar t, \bar x)$, we have:
\begin{equation*} 
\lim_{n \to \infty} d_{H}(t_n^\delta C(x_n), \bar \beta(t) C(\bar x)) = 0.
\end{equation*} 
Thus, Hausdorff continuity is established.

(iii) For any $x \in \widetilde C(t,u)$, there exists $(\theta_1,\cdots,\theta_m)^\top \in \Delta ^m$ such that
\begin{equation*} 
\begin{aligned}
\|x\|&=\left\|\beta(t)\sum_{i=1}^m\theta_i \nabla f_i(u)\right\|\le \beta(t) \max_{\theta \in \Delta^m}\left\|\sum_{i=1}^m \theta_i \nabla f_i(u)\right\|\\
&\le \beta(t) \left(\max_{\theta\in \Delta^m}\left\|\sum_{i=1}^m\theta_i (\nabla f_i(u)-\nabla f_i(0))\right\|+\max_{i\in[m]}\|\nabla f_i(0)\|\right)\\
&\le \beta(t) \cdot \left(L\|u\|+\max_{i\in[m]}\|\nabla f_i(0)\|\right).
\end{aligned}
\end{equation*} 
Taking $\varphi(t)=\beta(t) \cdot \max\left\{L, \max_{i\in[m]}\|\nabla f_i(0)\|\right\}$, we obtain the desired result.
\end{proof}
Based on the set-valued mapping $B$, we can define the following set-valued mapping:
\begin{equation}\label{eq:set-G}
\begin{aligned}
G: & [t_0, +\infty) \times \mathbb{R}^n \times \mathbb{R}^n \rightrightarrows \mathbb{R}^n \times \mathbb{R}^n, \\
& (t, u, v) \mapsto \{v\} \times \left( -\alpha(t)v - \mathop{\argmin}_{g \in B(t, u, v)} \langle g, -v \rangle \right).
\end{aligned}
\end{equation}
\begin{lemma}
Let the set-valued mapping $B$ be defined as in \cref{eq:set-B} and suppose \cref{assum:A-G-condition} holds. Then the set-valued mapping $I:(t,u,v)\mapsto \argmin_{g\in B(t,u,v)}\langle g,-v\rangle $ is upper semi-continuous.
\end{lemma}
\begin{proof}
According to \cite[Theorem 3B.5]{dontchev2009implicit}, the mapping $(t,u,v) \to \min_{g \in B(t,u,v)} \langle g, -v \rangle$ is continuous. We now proceed by contradiction. Suppose $I$ is not upper semi-continuous. Then there exists $(\bar{t}, \bar{u}, \bar{v})$ and some $\varepsilon_0 > 0$ such that for every $k > 0$,  
\begin{equation*} 
I\bigl(\mathcal{B}_{\frac{1}{k}}(\bar{t}, \bar{u}, \bar{v})\bigr) \not\subseteq I(\bar{t}, \bar{u}, \bar{v}) + \mathcal{B}_{\varepsilon_0}(0),
\end{equation*}   
so there exist sequences $\{g_k\}$ and $\{(t_k, u_k, v_k)\}$ satisfying $g_k \in I(t_k, u_k, v_k)$ with $(t_k, u_k, v_k) \in \mathcal{B}_{\frac{1}{k}}(\bar{t}, \bar{u}, \bar{v})$, and $d\bigl(g_k, I(\bar{t}, \bar{u}, \bar{v})\bigr) > \varepsilon_0$.

On the other hand, since $B$ is Hausdorff continuous, for the given $\varepsilon_0 > 0$ there exists $\delta > 0$ such that whenever $(t, u, v) \in \mathcal{B}_{\delta}(\bar{t}, \bar{u}, \bar{v})$,  
\begin{equation*} 
 d\bigl(a, B(\bar{t}, \bar{u}, \bar{v})\bigr) < d_H\bigl(B(t, u, v), B(\bar{t}, \bar{u}, \bar{v})\bigr) < \varepsilon_0.
\end{equation*}   
That is, $B\bigl(\mathcal{B}_{\delta}(\bar{t}, \bar{u}, \bar{v})\bigr) \subseteq B(\bar{t}, \bar{u}, \bar{v}) + \mathcal{B}_{\varepsilon_0}(0)$. Consequently, for sufficiently large $k$,  
\begin{equation*} 
g_k \in I\bigl(\mathcal{B}_{\frac{1}{k}}(\bar{t}, \bar{u}, \bar{v})\bigr) \subseteq B(\bar{t}, \bar{u}, \bar{v}) + \mathcal{B}_{\varepsilon_0}(0),
\end{equation*}   
meaning the sequence $\{g_k\}$ is bounded and therefore possesses a cluster point $\bar{g}$. Without loss of generality, assume $\lim_{k \to \infty} g_k = \bar{g}$. By continuity, we obtain  
\begin{equation*} 
\langle \bar{g}, -\bar{v} \rangle = \lim_{k \to \infty} \langle g_k, -v_k \rangle = \lim_{k \to \infty} \min_{g \in B(t_k, u_k, v_k)} \langle g, -v_k \rangle = \min_{g \in B(\bar{t}, \bar{u}, \bar{v})} \langle g, -\bar{v} \rangle.
\end{equation*}   
Thus, $\bar{g} \in I(\bar{t}, \bar{u}, \bar{v})$, which leads to the following contradiction:  
\begin{equation*} 
0 = d\bigl(\bar{g}, I(\bar{t}, \bar{u}, \bar{v})\bigr) = \lim_{k \to \infty} d\bigl(g_k, I(\bar{t}, \bar{u}, \bar{v})\bigr) > \varepsilon_0.
\end{equation*}   
This completes the proof.
\end{proof}
This set-valued mapping satisfies the following lemma:
\begin{lemma}\label{lem:solution-yesderday}
Let the set-valued mapping $B$ be defined as in \cref{eq:set-B} and suppose \cref{assum:A-G-condition} holds. Let $G$ be defined as in \cref{eq:set-G}. Then:
\begin{enumerate}[label=(\roman*)] 
\item[$\rm(i)$] For every $(t,u,v) \in [t_0, +\infty) \times \mathbb{R}^n \times \mathbb{R}^n$, the set $G(t,u,v) \subseteq \mathbb{R}^n \times \mathbb{R}^n$ is compact and convex;
\item[$\rm(ii)$] $G$ is upper semicontinuous;
\item[$\rm(iii)$] The mapping
\begin{equation*} 
\begin{aligned}
\phi: & [t_0, +\infty) \times \mathbb{R}^n \times \mathbb{R}^n \to \mathbb{R}^n \times \mathbb{R}^n,\  (t, u, v) \mapsto \mathop{\proj}_{G(t,u,v)}(0),
\end{aligned}
\end{equation*} 
is locally compact, i.e., for any point $x_0$, there exists a neighborhood $\mathcal{U}$ of $x_0$ such that $\phi(\mathcal{U})$ is compact.
\end{enumerate}
\end{lemma}
\begin{proof}
$\rm(i)$ For any $(t,u,v)$, according to \cref{assum:A-G-condition}, $B(t,u,v)$ is a compact convex set. Consequently, $\argmin_{g\in B(t,u,v)}\langle g,-v\rangle$ is a compact convex set. Since Cartesian products preserve convexity and compactness, $G(t,u,v)$ is compact and convex.

$\rm(ii)$ 
We now prove that for any $(\bar t,\bar u,\bar v)\in[t_0,+\infty) \times \mathbb{R}^n\times \mathbb{R}^n$ and any $\varepsilon >0$, there exists $\delta >0$ such that $G(\mathcal B_{\delta}(x_0)) \subseteq G(\bar t,\bar u,\bar v) + \mathcal B_{\varepsilon}(0)$.  

Since $I$ is upper semicontinuous, there exists $\delta_1$ such that for any $(t,u,v)$ satisfying $\|(t,u,v)-(\bar t,\bar u,\bar v)\|<\delta_1$ and any $g\in I(t,u,v)$, there exists $\bar g\in I(\bar t,\bar u,\bar v)$ with $\|g-\bar g\|<\frac{\varepsilon}{3}$.  

There exists $\delta_2$ such that for any $(t,u,v)$ satisfying $|t-\bar t|\le\|(t,u,v)-(\bar t,\bar u,\bar v)\|<\delta_2$, we have $|\alpha(t)-\alpha(\bar t)|\|v\|<\frac{\varepsilon}{3}$.  Let $M_{\delta_2}:=\sup_{t\in[\bar t-\delta_2,\bar t+\delta_2]} \alpha(t)$.

There exists $\delta_3$ such that for any $(t,u,v)$ satisfying $\|v-\bar v\|\le\|(t,u,v)-(\bar t,\bar u,\bar v)\|<\delta_3$, we have $\|v-\bar v\|<\frac{\varepsilon}{3(1+M_{\delta_2})}$.  

Taking $\delta = \min\{\delta_1,\delta_2,\delta_3\}$, for any $(t,u,v)$ satisfying $\|(t,u,v)-(\bar t,\bar u,\bar v)\|<\delta$ and any $(x,y)=(v,-\alpha(t)v-g)\in G(t,u,v)$, there exists $(\bar x,\bar y)=(\bar v,-\alpha(\bar t)\bar v-\bar g)\in G(\bar t,\bar u,\bar v)$ such that
\begin{equation*} 
\begin{aligned}
\|(x,y)-(\bar x,\bar y)\|&\le \|v-\bar v\|+\left\|-\alpha(t)v-g+\alpha(\bar t)\bar v+\bar g\right\|\\
&\le \|v-\bar v\| + \left\|-\alpha(t)v+\alpha(\bar t)\bar v\right\| + \left\|-g+\bar g\right\|\\
&\le \|v-\bar v\| + |\alpha(t)|\|v-\bar v\| + |\alpha(t)-\alpha(\bar t)|\|\bar v\| + \|g-\bar g\|\\
&\le (1+\alpha(t))\|v-\bar v\| + |\alpha(t)-\alpha(\bar t)|\|\bar v\| + \|g-\bar g\|\\
&< \varepsilon.
\end{aligned}
\end{equation*} 
This proves the upper semicontinuity.

$\rm(iii)$ The result follows from the upper semicontinuity.
\end{proof}

Based on \cref{eq:set-G}, we can define the differential inclusion:
\begin{equation}\label{eq:DI}
\left\{\begin{aligned}
&(\dot u(t),\dot v(t)) \in G(t,u(t),v(t)),~~\text{ for }t>t_0,\\
&(u(t_0),v(t_0))=(u_0,v_0).
\end{aligned}
\right. \tag{DI}
\end{equation}
\begin{theorem}[{\cite[Theorem 3.2]{sonntag2024fastSIAMonOptimization,aubin2009differential}}]\label{thm:exist}
    Let $\mathcal{X}$ be a real Hilbert space, and let $\Omega \subset \mathbb{R} \times \mathcal{X}$ be an open subset containing $(t_0, x_0)$. Let $G$ be an upper semicontinuous map from $\Omega$ into the nonempty closed convex subsets of $\mathcal{X}$. We assume that $(t, x) \mapsto \proj_{G(t,x)}(0)$ is locally compact. Then, there exists $T > t_0$ and an absolutely continuous function $x$ defined on $[t_0, T]$, which is a solution to the differential inclusion
    \begin{equation*}
    \dot{x}(t) \in G(t, x(t)), \quad x(t_0) = x_0.
    \end{equation*}
\end{theorem}
\begin{theorem}\label{thm:existence-theorem}
For any $(u_0,v_0)\in \mathbb{R}^n \times \mathbb{R}^n$, there exists $T>0$ and an absolutely continuous function $(u,v)$ defined on $[t_0,T]$ that is a solution of \cref{eq:DI} on $[t_0,T]$.
\end{theorem}
\begin{proof}
According to \cref{lem:solution-yesderday,thm:exist}, we obtain the conclusion.
\end{proof}

\begin{theorem}\label{thm:global-solution}
{Let $B$ be the set-valued mapping defined in \cref{eq:set-B} and satisfying \cref{assum:A-G-condition}. For any $(u_0, v_0) \in \mathbb{R}^n \times \mathbb{R}^n$, there exists a smooth function $x: [t_0, +\infty) \to \mathbb{R}^n$ that is a solution in the sense of \cref{def:soultion-general}.}
\end{theorem}
\begin{proof}
Define the set
\begin{equation*} 
\begin{aligned}
\mathfrak{S} := \{(u,v,T): & ~T \in (t_0, +\infty] \text{ and } (u,v): [t_0, T) \to \mathbb{R}^n \times \mathbb{R}^n \text{ is absolutely continuous on every} \\
&\text{compact interval contained in $[t_0, T)$ and is a solution of \cref{eq:DI} on $[t_0, T)$} \}.
\end{aligned}
\end{equation*}
By \cref{thm:existence-theorem}, $\mathfrak{S}$ is nonempty. We define a partial order $\preccurlyeq$ on $\mathfrak{S}$ as follows: for $(u_1,v_1,T_1), (u_2,v_2,T_2) \in \mathfrak{S}$,
\begin{equation*} 
(u_1,v_1,T_1) \preccurlyeq (u_2,v_2,T_2) \quad \Longleftrightarrow \quad T_1 \le T_2 \text{ and } (u_1(t),v_1(t)) = (u_2(t),v_2(t)) \text{ for all } t \in [t_0, T_1).
\end{equation*} 
Following the approach in \cite{sonntag2024fastSIAMonOptimization}, by Zorn's Lemma, there exists a maximal element $(u,v,T)$ in $\mathfrak{S}$. If $T = +\infty$, the proof is complete.

Assume $T < +\infty$. Define the function on $[t_0, T)$:
\begin{equation*} 
h(t) := \|(u(t),v(t)) - (u(t_0),v(t_0))\|.
\end{equation*} 
Then, for almost all $t \in [t_0, T)$,
\begin{equation*} 
\frac{d}{dt}\left(\frac12 h^2(t)\right) \le \|(\dot u(t),\dot v(t))\| h(t).
\end{equation*} 
Since $(\dot u, \dot v) \in G(t,u,v)$ on $[t_0, T)$, we have
\begin{equation*} 
\begin{aligned}
\|(\dot u(t),\dot v(t))\| &= \left\|\left(v, -\alpha(t)v - g\right)\right\| \le \|v\| + \alpha(t)\|v\| + \|g\|\\& \le \left(1 + \sup_{t\in [t_0,T]} \alpha(t)\right) \|v\| + \varphi(t)(1 + \|(u,v)\|) \le c(1 + \|(u,v)\|),
\end{aligned}
\end{equation*} 
where $c = \sqrt{2} \max\left\{\left(1 + \sup_{t\in[t_0,T]} \alpha(t)\right), \sup_{t \in [t_0, T)} \varphi(t)\right\}$. By the triangle inequality,
\begin{equation}\label{eq:diff-bounded}
\|(\dot u(t),\dot v(t))\| \le \widetilde{c}(1 + \|(u(t),v(t)) - (u(t_0),v(t_0))\|),
\end{equation}
where $\widetilde{c} = c(1 + \|(u(t_0),v(t_0))\|)$. Therefore, for almost all $t \in [t_0, T)$,
\begin{equation*} 
\frac{d}{dt}\left(\frac12 h^2(t)\right) \le \widetilde{c}(1 + h(t)) h(t).
\end{equation*}
Integrating over $[t_0, t)$ yields
\begin{equation*} 
\int_{t_0}^t \frac{1}{2} \frac{d}{ds} h^2(s)  ds \le \int_{t_0}^t \widetilde{c}(1 + h(s)) h(s)  ds,
\end{equation*} 
so for all $t \in [t_0, T)$,
\begin{equation*} 
\frac12 h(t)^2 \le \int_{t_0}^t \widetilde{c}(1 + h(s)) h(s)  ds.
\end{equation*}
Applying \cite[Lemma A.2]{attouch2015multiibjective}, for all $t \in [t_0, T)$,
\begin{equation*}
\begin{aligned}
h(t) &\le \int_{t_0}^t \widetilde{c}(1 + h(s)) ds \le \widetilde{c} T + \int_{t_0}^t \widetilde{c} h(s)  ds.
\end{aligned}
\end{equation*}
Applying \cite[Lemma A.1]{attouch2015multiibjective}, for all $t \in [t_0, T)$,
\begin{equation*} 
h(t) \le \widetilde{c} T \exp(\widetilde{c} T).
\end{equation*} 
Thus, $h$ is bounded on $[t_0, T)$. Therefore, $u$ and $v$ are also bounded on $[t_0, T)$. By \cref{eq:diff-bounded}, $\dot u$ and $\dot v$ are essentially bounded. Since $u$ and $v$ are absolutely continuous, for all $t \in [t_0, T)$,
\begin{equation*} 
\begin{aligned}
u(t) &= u_0 + \int_{t_0}^t \dot u(s)  ds \in \mathbb{R}^n, &
v(t) &= v_0 + \int_{t_0}^t \dot v(s)  ds \in \mathbb{R}^n.
\end{aligned}
\end{equation*} 
Hence, for any $t_1, t_2$,
\begin{equation*} 
\begin{aligned}
\|u(t_1) - u(t_2)\| &\le \operatorname{esssup}_{s \in [t_1,t_2]} \|\dot u(s)\| \cdot |t_1 - t_2|, \\
\|v(t_1) - v(t_2)\| &\le \operatorname{esssup}_{s \in [t_1,t_2]} \|\dot v(s)\| \cdot |t_1 - t_2|.
\end{aligned}
\end{equation*}
Therefore, the limits $v_T := \lim_{t \to T^-} v(t) \in \mathbb{R}^n$ and $u_T := \lim_{t \to T^-} u(t) \in \mathbb{R}^n$ exist. Now consider the differential inclusion:
\begin{equation*} 
\left\{\begin{aligned}
&(\dot u(t), \dot v(t)) \in G(t, u(t), v(t)) \quad \text{for } t > T, \\
&(u(T), v(T)) = (u_T, v_T),
\end{aligned}\right.
\end{equation*}
By \cref{thm:existence-theorem}, there exists $\delta > 0$ and a solution $(\tilde u, \tilde v): [T, T+\delta] \to \mathbb{R}^n \times \mathbb{R}^n$ that is absolutely continuous on any bounded closed interval in $[T, T+\delta]$. Define
\begin{equation*} 
(u^*, v^*): [t_0, T+\delta) \to \mathbb{R}^n \times \mathbb{R}^n, \quad t \mapsto \left\{\begin{aligned}
&(u(t), v(t)),& &\text{for } t \in [t_0, T), \\
&(\tilde u(t), \tilde v(t)),& &\text{for } t \in [T, T+\delta).
\end{aligned}\right.
\end{equation*}
Then $(u^*, v^*, T+\delta) \in \mathfrak{S}$ and $(u,v,T) \preccurlyeq (u^*,v^*,T+\delta)$ with $(u,v,T) \neq (u^*,v^*,T+\delta)$. This contradicts the maximality of $(u,v,T)$ in $\mathfrak{S}$.
\end{proof}
\begin{lemma}[{\cite[Lemma 3.8]{sonntag2024fastSIAMonOptimization}}]\label{lem:son}
    Let $C\subset \R^n$ a convex and compact set, and $\eta \in \R^n$ a fixed vector. Then $\xi\in \eta- \argmin_{\mu\in C}\langle \mu,\eta\rangle  $ if and only if $\eta =\proj_{C+\xi}(0)$
\end{lemma}
\begin{theorem}
Suppose the set-valued mapping $B: (t_{\min}, +\infty) \times \mathbb{R}^n \times \mathbb{R}^n \rightrightarrows \mathbb{R}^n$ satisfies \cref{assum:A-G-condition}. Let $t_0 > 0$ and $x_0, v_0 \in \mathbb{R}^n$. If $(u, v): [t_0, \infty) \to \mathbb{R}^n \times \mathbb{R}^n$ is a solution of \cref{eq:DI} with $(u(t_0), v(t_0)) = (u_0, v_0)$, then $x(t) := u(t)$ satisfies the following equation:
$$
\alpha(t) \dot x(t) + \proj_{B(t, x(t), \dot x(t)) + \ddot x(t)}(0) = 0,
$$
for almost all $t \in [t_0, +\infty)$, with $x(t_0) = x_0$ and $\dot x(t_0) = v_0$, i.e. satisfies \cref{eq:CP}. 
\end{theorem}
\begin{proof}
By \cref{thm:global-solution}, the conclusion follows.
\end{proof}
\section*{Acknowledgements} 
We are grateful to the referees for their constructive comments and suggestions, which have greatly improved the presentation of this paper.

\end{document}